\documentclass{article}

\usepackage{arxiv}
\usepackage[title]{appendix}
\usepackage{graphicx}
\usepackage[pdfborder={0 0 0}]{hyperref}
\usepackage{bm}
\usepackage{cite}
\usepackage{comment}
\usepackage{graphics,fancyhdr,graphicx,subfigure}
\usepackage{subfigure}
\usepackage{cancel}
\usepackage{amsmath,amsfonts,latexsym,amssymb}
\usepackage{lineno}
\usepackage[ruled,linesnumbered]{algorithm2e}
\usepackage{multirow}
\usepackage{listings}
\usepackage[table]{xcolor}
\usepackage{lscape}
\usepackage{tabularx}
\DeclareMathOperator*{\argmin}{\arg\!\min}

\def\equationautorefname~#1\null{%
  Eq.~(#1)\null
  }
\def\subfigureautorefname~#1\null{%
  Fig.~#1\null
}

\definecolor{listinggray}{gray}{0.9}
\definecolor{lbcolor}{rgb}{0.9,0.9,0.9}
\definecolor{Darkgreen}{RGB}{0,100,0}

\usepackage{array}
\newcolumntype{L}[1]{>{\raggedright\let\newline\\\arraybackslash\hspace{0pt}}m{#1}}
\newcolumntype{C}[1]{>{\centering\let\newline\\\arraybackslash\hspace{0pt}}m{#1}}
\newcolumntype{R}[1]{>{\raggedleft\let\newline\\\arraybackslash\hspace{0pt}}m{#1}}

\title{Threshold shift method for reliability-based design optimization}

\author{
  Somdatta~Goswami\thanks{These authors have equal contributions} \\
  Institute of Structural Mechanics\\ Bauhaus-Universit{\"a}t \\
  Weimar, 99423-Weimar, Germany\\
  \texttt{somdatta.goswami@uni-weimar.de} \\
   \And
 Souvik~Chakraborty${^*}$ \\
  Department of Aerospace and Mechanical Engineering\\
  University of Notre Dame\\
  Notre Dame, IN - 46556, U.S.A. \\
  \texttt{csouvik41@gmail.com; schakrab@nd.edu} \\
   \And
 Rajib~Chowdhury \\
  Department of Civil Engineering\\
  Indian Institute of Technology Roorkee\\
  Roorkee, India \\
  \texttt{rajibfce@iitr.ac.in} \\
   \And
 Timon~Rabczuk \\
  Department of Computer Engineering\\
  College of Computer and Information Sciences\\
  King Saud University, Riyadh, Saudi Arabia\\
  \texttt{timon.rabczuk@uni-weimar.de}
}

\begin{document}
\maketitle

\begin{abstract}
We present a novel approach, referred to as the ‘threshold shift method’ (TSM), for reliability based design optimization (RBDO). The proposed approach is similar in spirit with the sequential optimization and reliability analysis (SORA) method where the RBDO problem is decoupled into an optimization and a reliability analysis problem. However, unlike SORA that utilizes {\it shift-vector} to shift the design variables within a constraint (independently), in TSM we propose to shift the threshold of the constraints. We argue that modifying a constraint, either by shifting the design variables (SORA) or by shifting the threshold of the constraints (TSM), influences the other constraints of the system. Therefore, we propose to determine the thresholds for all the constraints by solving a single optimization problem. Additionally, the proposed TSM is equipped with an active-constraint determination scheme. 
To make the method scalable, a practical algorithm for TSM that utilizes two surrogate models is proposed. Unlike the conventional RBDO methods, the proposed approach has the ability to handle highly non-linear probabilistic constraints. The performance of the proposed approach is examined on six benchmark problems selected from the literature. The proposed approach yields excellent results outperforming other popular methods in literature. As for the computational efficiency, the proposed approach is found to be highly efficient, indicating it's future application to other real-life problems.
\end{abstract}

\section{Introduction}
\label{sec:intro}
Uncertainty is an eternal companion of all structural systems.
In a system, uncertainty can originate from
both, uncertain parameters within the model as well as the limitations in the model itself.
These uncertainties, if ignored, may result in to catastrophic failures.
Therefore, it  is  necessary  to  consider  the  effect  of  uncertainty  in  an optimization process.
In literature, there exist two methodologies to incorporate uncertainty into the framework of optimization, namely the robust design optimization (RDO) \cite{Chakraborty2017a, Chatterjee2019a,Youn2009} and the reliability-based design optimization (RBDO) \cite{Dubourg2013,Dubourg2011,Dubourg2011a, Moon2018c, Meng2018c, Bichon2008}. In this study, we focus on RBDO. 

The primary bottleneck of RBDO resides in the fact that one has to solve a reliability analysis problem within each iteration of a RBDO problem.
Quite naturally, this becomes computationally expensive.
Moreover, computing the gradients of the probabilistic constraints is also challenging.
To address these issues, different RBDO formulations have been proposed by the researchers.
Broadly speaking, the RBDO methods available in literature can be categorized into two groups, namely the two-level approaches \cite{Youn2003,Youn2004} and the single-loop approaches \cite{Liang2008,Keshtegar2018a,Thanedar1992, Meng2018b, Meng2018a, meng2019, Picheny2010}.

As the name suggests, the two-level approaches involve two stages of analyzing the problem; one corresponding to the reliability analysis and the other corresponding to the design optimization.
One popular two-level approach is the reliability index approach (RIA)\cite{Shi2016,Yu1998,Enevoldsen1994,Yu1997}.
In RIA, we represent the probabilistic constraints in terms of the reliability indices and solve a coupled problem.
Within each iteration of a RBDO problem, the reliability analysis problem is generally solved by using the well-known first-order reliability analysis method (FORM) \cite{Lopez2012,Hu2015first,Du2012first,Tichy1994first,Choi2012uncertainty}.
Unfortunately, RIA is often unstable and the rate of convergence is rather slow.
Additionally, results obtained are often not optimal, 
specifically for problems involving
highly non-linear probabilistic constraints.

In order to address the issues associated with RIA, one can modify the method by introducing a surrogate model. 
The primary idea is to replace the computationally expensive limit state function with a computationally efficient surrogate models.
Popular surrogate models available in literature include polynomial chaos expansion \cite{Xiu2002athe, Sudret2008global}, Kriging \cite{Bilionis2013multi,Biswas2018saturation,Biswas2016kriging,Mukhopadhyay2016a}, analysis of variance decomposition \cite{Chakraborty2017polynomial, Chakraborty2016hybrid, Chakraborty2017towards, Chakraborty2017an}, neural network \cite{HosniElhewy2006, Deng2005} and support vector machine \cite{Zhao2014, Guo2009}.
Use of surrogate models within the RBDO framework can be found in \cite{Dubourg2011,Dubourg2011a,Rahman2010reliability,Rahman2008design, leger2017reliability}.

Another popular two-level approach for RBDO is the performance measure approach (PMA) \cite{Shi2016,Tu1999,Youn2005}.
In PMA, the probabilistic constraints are transformed to performance constraints \cite{Youn2004}.
Unlike RIA that solves a simple optimization problem under complicated probabilistic constraints, PMA solves a complicated optimization problem subjected to simple constraints.
As a consequence, PMA is somewhat more stable compared to RIA.
Unfortunately, PMA is also computationally expensive as one have to solve several inverse reliability analysis problems. Improvements to PMA include chaos control method \cite{Yang2009}, hybrid chaos control performance measure approach \cite{Meng2015}, self-adaptive modified chaos control \cite{Keshtegar2017}, enhanced chaos control method \cite{Hao2017, hao2017efficient} and  dynamical accelerated chaos control method \cite{Keshtegar2018}.

As an improvement of PMA, the sequential optimization and reliability analysis (SORA) method was proposed in \cite{Du2004}. 
In SORA, the RBDO problem is decoupled such that, within each iteration, one has to sequentially solve a deterministic optimization problem and a reliability analysis problem.
To that end, the concept of {\it shift-vector} is introduced.
Unfortunately, the dimensionality of the shift-vector increases with increase in number of design variables and hence, the application of SORA is limited to academic problems. 
Modifications to conventional SORA such as the one presented in \cite{hao2019new} can also be found in the literature.

In recent times, the single loop approaches have become quite popular among researchers. 
In this approach one attempts to convert the double-loop procedure into a single one. 
To that end, one option is to replace the probabilistic constraints by optimality conditions \cite{Liang2008}. This is achieved by converting a probabilistic constraint to a deterministic constraint \cite{Li2013}. 
Compared to the two-level approaches mentioned earlier, the single loop approaches are computationally efficient.
However, converting the nested loop into a single loop often yields erroneous results.
Moreover, single loop approaches can also lead to numerical instability \cite{chen1997,Nguyen2010,Li2013,Liang} as it does not satisfy the Karush-Kuhn-Tucker necessary condition. 

To address the issues associated with the available methods, we propose a novel method for RBDO.
The proposed approach, referred to as threshold shift method (TSM), is an improvement over SORA, proposed in \cite{Torii2016a}. However, unlike SORA, TSM does not require a shift-vector for the design variables; instead, the thresholds associated with the constraints are shifted.
We argue that modifying  a  constraint, either by shifting the design variables \cite{Torii2016a} or by shifting the threshold of the constraints (TSM), influence the other constraints in the system.  
Therefore, we propose to determine the threshold of all the constraints involved in the system by solving a single optimization problem. 
To make TSM scalable to problems having large number of constraints, an active-constraint determination scheme is proposed and integrated with TSM. Lastly, to increase the efficiency of optimization and reliability analysis, two surrogate models have been utilized within the TSM framework. It is to be noted that unlike FORM based RBDO approaches, TSM can solve problems with highly non-linear probabilistic constraints.

The remainder of the paper is organized as follows. In \autoref{sec:RBDO} we provide a brief description of the fundamental concepts of RBDO and an inception to SORA. \autoref{sec:sora_alternate} discusses an alternative view of SORA, which forms the basis of the present work. The details about the proposed approach is presented in \autoref{sec:methodology}. Numerical examples illustrating the performance of the proposed approach are presented in \autoref{sec:numericals}. Finally, \autoref{sec:conclusion} presents the concluding remarks.

\section{Reliability based design optimization}
\label{sec:RBDO}
Mathematically, a RBDO problem can be defined as:
\begin{equation}\label{eq:RBDO_model}
  \begin{split}
        & \text{minimize        } c(\bm{d})\\
        & \text{s.t. } \\
        & \mathbb{P}(G_i \left(\bm{d},\bm{X})< 0\right)) -\Phi(-\beta_{\text{Tar}})\leq 0 \text{   ,   } i = 1,\ldots,n_p,\\
        & \bm{d}^{L}\leq \bm{d}\leq \bm{d}^{U},\,
        \end{split}
\end{equation}
where $c$ represents the objective/cost function, $\bm{d} = \left(d_1, \ldots, d_m \right) \in \mathbb{R}^m$ and $\bm{X} = \left(X_1, \ldots, X_n \right) \in \mathbb{R}^n$ represents the design variables and the random variables, respectively and $n_p$ is the number of probabilistic constraints. $G\left(\bm{d}, \bm{X} \right)$ in \autoref{eq:RBDO_model} denotes the limit state function such that $G\left(\bm{d}, \bm{X} \right) < 0$ is the failure domain. $\mathbb{P} \left( \bullet \right)$ represents the probability measure and $\beta_{\text{Tar}}$ represents the target reliability index. Often, we consider $\bm{d} = \bm{\mu_{X}}$, where $\bm{\mu_X}$ represent the mean of the random variables. It is to be noted that a RBDO problem may also have some non-linear deterministic constraints.

\noindent The probability of failure in \autoref{eq:RBDO_model} is computed as
\begin{equation}\label{eq:pf_compute}
\begin{split}
    \mathbb{P}\left(G_i\left( \bm{X}, \bm{d} \right) \le 0\right) &= \int_{\Omega_f}{f_{\bm X}\left(\bm{x} \right)} d\bm{x}\\
    &= \int_{\Omega}{\mathcal{I}\left(G_i\left( \bm{X}, \bm{d} \right) < 0 \right)f_{\bm X}\left(\bm{x}\right)}d\bm{x},
\end{split}
\end{equation}
where $\Omega_f$ and $\Omega$ indicate the failure domain and the problem domain, respectively and $\mathcal{I}$ is an indicator function. $f_{\bm X}$ denotes the probability density function of the random variables.
Note that the integration in \autoref{eq:pf_compute} is high-dimensional and approximate methods have to be used for computing the probability of failure (i.e., structural reliability analysis).
Unfortunately, the simulation based approaches (e.g.Monte Carlo simulation (MCS) and Importance simulation) are computationally expensive and do not provide estimate of the sensitivities. As a result these are not that suitable for RBDO problems.
Use of first and second-order reliability methods in RBDO can be found in literature. In fact, both these methods have been used within the framework of RIA and PMA.
Unfortunately both RIA and PMA have several limitations as highlighted in \autoref{sec:intro}. 

\begin{figure}
    \centering
    \includegraphics[width= 0.45\textwidth]{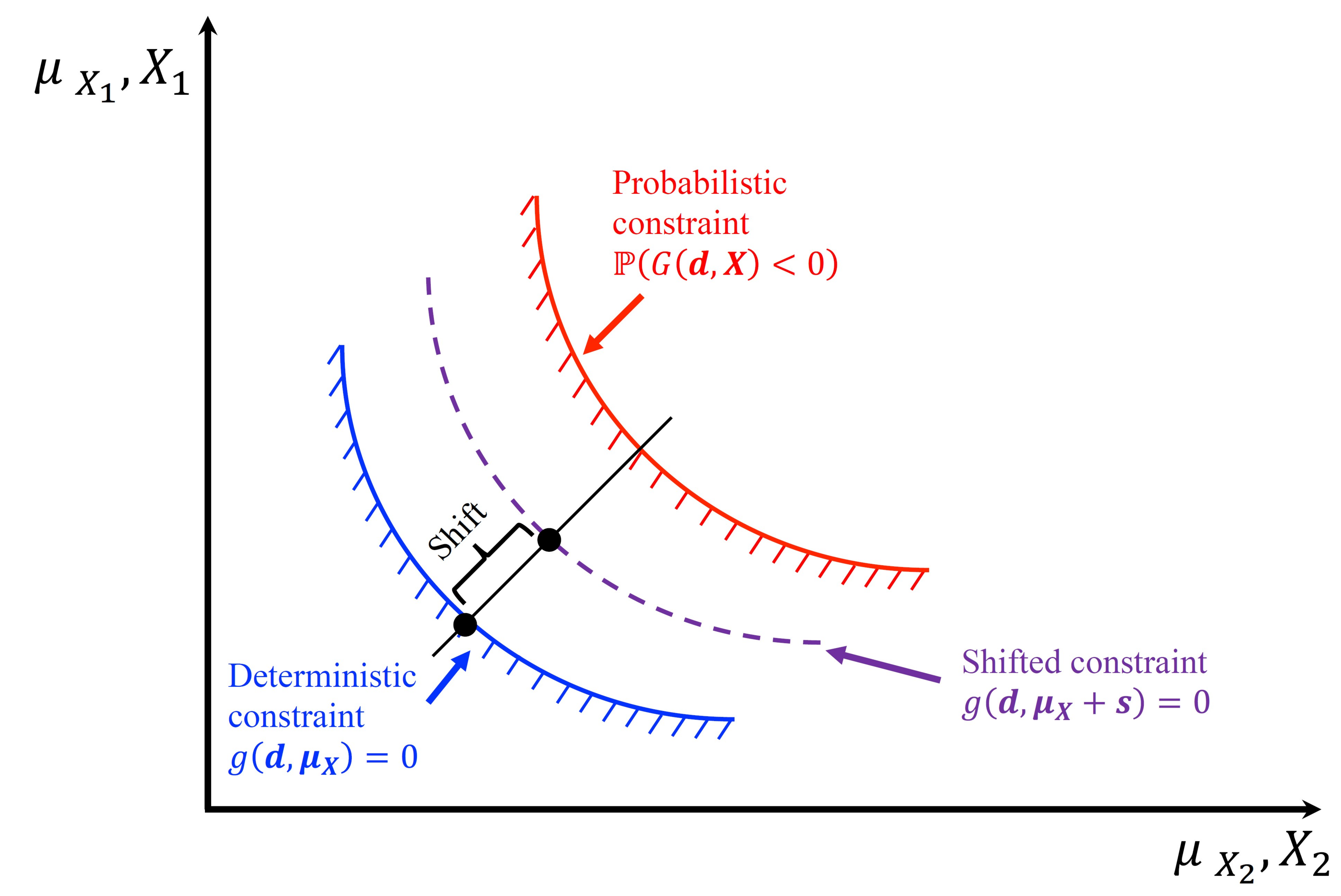}
    \caption{Schematic representation of the SORA method. The blue line represents the deterministic constraint and the red line represents the probabilistic constraint. The violet (dashed) line represents the shifted constraint. At convergence, the shifted constraint will overlap or be extremely close to the probabilistic constraint.}
    \label{fig:SORA}
\end{figure}

An interesting alternative to RIA and PMA is the SORA.
The primary idea of the SORA approach is to shift the boundary of the violated equivalent deterministic constraint towards a feasible direction.
The feasible direction is determined based on the information obtained from the reliability analysis step.
A schematic representation of the SORA approach is shown in~\autoref{fig:SORA}. The optimization converges when the difference between the shifted equivalent-deterministic constraint and the actual probabilistic constraint is within prescribed tolerance limits.
The primary advantage of SORA resides in the fact that instead of solving a probabilistic optimization problem, we solve a series of sequential deterministic optimization and reliability analysis problems. 
As a consequence, well-developed optimization and reliability analysis tools can be used for this purpose. Moreover, the computational cost associated with SORA is significantly less compared to the original problem described in \autoref{eq:RBDO_model}.
Further details on SORA and its implementation can be found in \cite{Du2004}.

Despite the advantages discussed above, SORA is susceptible to all the limitations associated with FORM and hence, may not yield accurate results for problems with highly non-linear probabilistic constraints.
The objective of this work is to present a novel method for solving the RBDO problem stated in \autoref{eq:RBDO_model}.

\section{SORA - an alternative view}
\label{sec:sora_alternate}
Recently, an alternative view of SORA has been presented.
This {\it alternate view} stems from the observation that
a given shift-vector $\bm s$ will result in a fixed optimum design $\bm d$.
As a consequence, it is possible to build a mapping from $\bm s$ to $\bm d$
\begin{equation}\label{eq:map1}
    O: \bm s \rightarrow \bm d.
\end{equation}
On the other hand, a given design vector $\bm d$ results in a fixed probability of failure, $P_f = \mathbb{P}\left( G \left( \bm X, \bm d \right) < 0 \right)$. Hence, a mapping from $\bm d$ to $P_f$ can also be built.
\begin{equation}\label{eq:map2}
    R: \bm d \rightarrow P_f.
\end{equation}
The mapping in \autoref{eq:map1} is based on the optimization problem defined in SORA \cite{Torii2016a}, whereas the mapping in \autoref{eq:map2} is based on a reliability analysis problem for a constant $\bm d$.
Combining Eqs. (\ref{eq:map1}) and (\ref{eq:map2}), a composite functional relation can be built as
\begin{equation}\label{eq:map3}
    J = R  \circ O : \bm s \rightarrow P_f,
\end{equation}
where $J$ represents the mapping from $\bm s$ to $P_f$.
Alternatively, \autoref{eq:map3} can be viewed as
\begin{equation}
    J : \bm s \rightarrow \beta,
\end{equation}
With this, the goal is to determine $\bm s^*$ such that $P_f = P_t$, where $P_t$ is the target probability of failure. 
When this is accomplished, the optimum design vector $\bm d^*$ can be obtained by solving an optimization problem similar to SORA.
Mathematically this can be represented as:
\begin{equation}\label{eq:sora_alt}
\begin{split}
    & \bm s^* = \argmin_{\bm s} \left[ \beta \left(\bm s \right) - \beta_{\text{Tar}} \right]^2, \\
    \text{s.t. }\;\;\;\; &\beta_{\text{Tar}} - \beta\left(\bm s \right) \le 0.
\end{split}
\end{equation}
This is an alternative view of the conventional SORA algorithm.
For further details on the alternative SORA approach, interested readers may refer to \cite{Torii2016a}.

\begin{figure}
    \centering
    \subfigure[with design variable $\bm d$]{
    \includegraphics[width=0.4\textwidth]{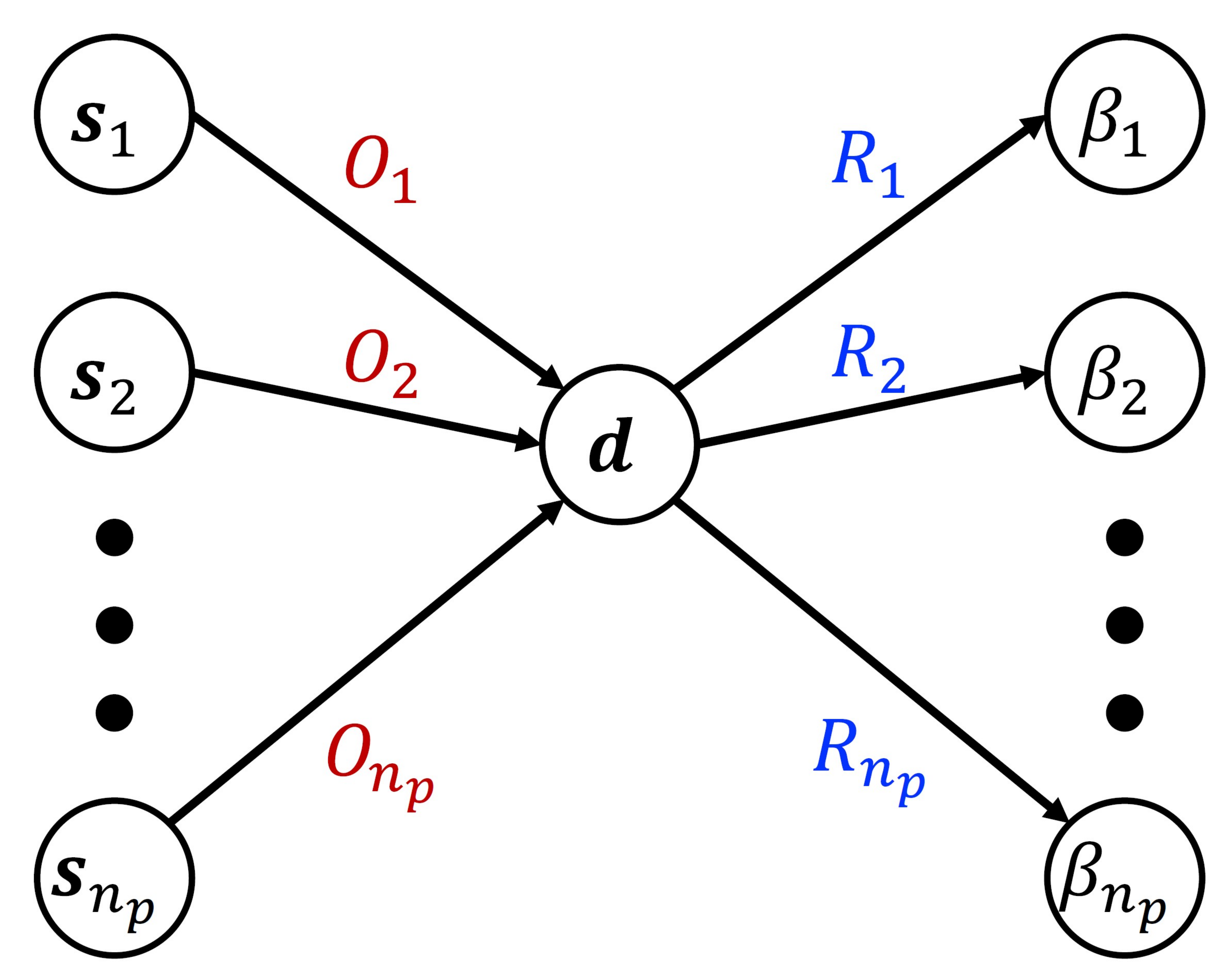}}
    \subfigure[direct mapping]{
    \includegraphics[width=0.4\textwidth]{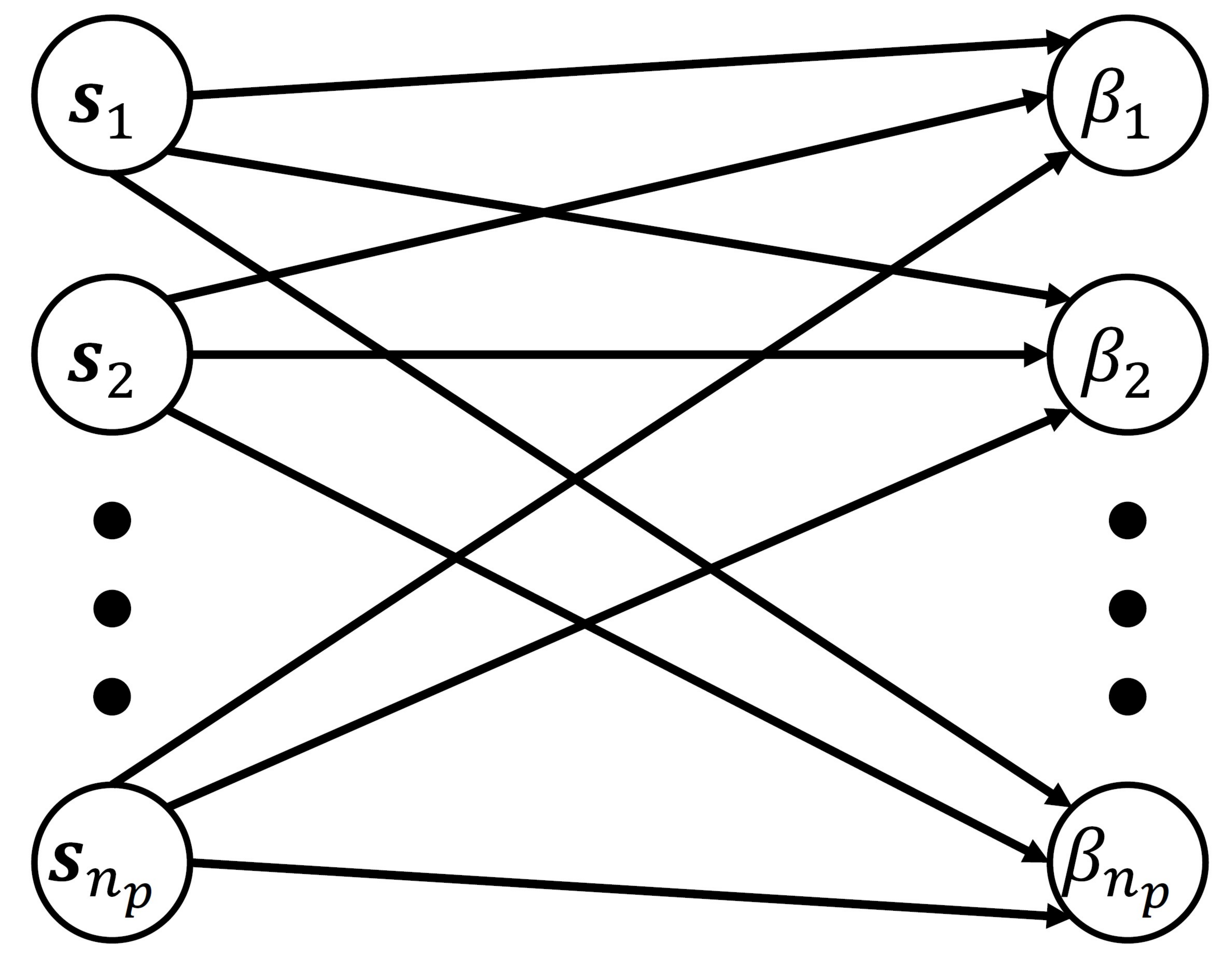}
    }
    \caption{Graphical models depicting the relations between the shift-vectors and the reliability indices. shift-vectors influence the reliability indices through the design variables (part a). In case of direct mapping (as in SORA), the design variable needs to be integrated out of the system. Under that scenario, the system becomes coupled (as shown in part b)}
    \label{fig:pgm_shift_vector}
\end{figure}

It was argued that for problems having several probabilistic constraints, one needs to solve~\autoref{eq:sora_alt} multiple times - one iteration per constraint \cite{Torii2016a}. 
It is to be noted that the above statement is an approximation.
The shift-vectors influences the reliability indices through the design variables, which in turn depends on all the shift-vectors.
A graphical representation of the same is shown in \autoref{fig:pgm_shift_vector}.

This alternate view also has certain associated challenges.
First, treating shift-vectors corresponding to probabilistic constraints independently is an approximation that can result in non-optimal solutions. Second, for problems with large number of random variables, the shift-vector will be of high-dimensional. This creates two issues:
    \begin{itemize}
        \item The optimization problem for computing the shift-vector (\autoref{eq:sora_alt}) becomes high dimensional and hence, computationally expensive.
        \item The mapping from $\bm s$ to $\beta$ also becomes computationally expensive.
    \end{itemize}
To address these issues we present a novel method, referred to as TSM. TSM can be viewed as a modification to the alternative view of SORA, proposed in \cite{Torii2016a}.
However, unlike the method proposed in \cite{Torii2016a}, we do not use shift vectors. Instead, the thresholds of the constraints are shifted. The details about the proposed TSM are furnished in the next section.

\section{Threshold shift method}
\label{sec:methodology}
In this section, the details of the proposed {\it threshold shift method} are presented along with its constituents. Towards the end of this section, we propose a dual-surrogate based approach for implementing TSM in practice.
\subsection{Proposed approach}
Consider $\bm{X} = \left(X_1, \ldots, X_n \right) \in \mathbb{R}^n$ to be the random variables associated with a system. 
Without any loss in generality, we also consider $\bm d = \bm {\mu_X}$. 
As discussed previously, the shift-vector for this problem would be an $n$-dimensional vector.
Proceeding a step further, if we assume that the underlying RBDO problem is having $n_p$ probabilistic constraints, one needs to solve $n_p$ optimization problems per iteration where each optimization problem involves $n$ design variables.
In cases where the number of random variables $n$ is large, the SORA algorithm will become computationally expensive.
Even for low-dimensional problems, the assumption that the shift-vectors can be computed independently may lead to erroneous solutions.

To address the above mentioned issues, TSM is proposed.
In TSM, instead of using a shift-vector for the design variables, we propose to shift the thresholds of the constraints.
With this, the optimization problem that previously involved $n$ design variables, now involves only one variable.
Formally in TSM, we solve an optimization problem of the form
\begin{equation}\label{eq:tsm1}
    \begin{split}
        &\bm d^* = \argmin _{\bm d}c\left(\bm d \right), \\
        \text{s.t. }\;\;\;\; & g\left(\bm d, \bm{\mu _X} \right) \le c_s \\
        & \bm d^L \le \bm d \le \bm d^U.
    \end{split}
\end{equation}
Note that \autoref{eq:tsm1} is similar to that of SORA, except that in place of the shift-vector $\bm s$, we introduce a threshold $c_s$.
The threshold $c_s$ in \autoref{eq:tsm1} is computed as
\begin{equation}\label{eq:tsm2}
\begin{split}
    & c_s^* = \argmin_{c_s} \left[\beta_{\text{Tar}} - \beta \left( c_s\right) \right]^2,\\
    \text{s.t. }\;\;\;\; &\beta_{\text{Tar}} - \beta\left(c_s\right) \le 0.
\end{split}
\end{equation}
Eqs. (\ref{eq:tsm1}) and (\ref{eq:tsm2}) are the fundamental equations of the proposed TSM.
Having said that it is to be noted that TSM, as presented {\it via} Eqs. (\ref{eq:tsm1}) and (\ref{eq:tsm2}), is applicable only to systems with single-probabilistic constraint.
It is possible to follow a similar procedure as proposed in SORA and extend this to problems with multiple probabilistic constraints.
However, as already stated, such an approach is an approximation and may yield erroneous results.
Therefore, in this work, we adopt an alternative path.
A coupled-optimization problem is solved to compute the thresholds corresponding to all the constraints simultaneously.
Mathematically this can be represented by modifying Eqs. (\ref{eq:tsm1}) and (\ref{eq:tsm2}) as follows:
\begin{equation}\label{eq:tsm_mc1}
    \begin{split}
        & \bm d^* = \argmin c\left(\bm d \right), \\
        \text{s.t. }\;\;\;\; & g_i\left(\bm d, \bm {\mu_X} \right) \le c_s^i, i=1,\ldots, n_p,\\
        & \bm d^L \le \bm d \le  \bm d^U,
    \end{split}
\end{equation}
\begin{equation}\label{eq:tsm_mc2}
    \begin{split}
        & \bm c_s^* = \argmin_{\bm c_s} \sum_{i=1}^{n_p}\left[ \beta_{\text{Tar},i} - \beta _ {i}\left(\bm c_s \right) \right]^2, \\
        \text{s.t. }\;\;\;\; & \beta_{\text{Tar},i} - \beta _ {i}\left(\bm c_s \right) \le 0, i=1,\ldots, n_p, \\
    \end{split}
\end{equation}
where $\bm c_s = \left(c_s^1,\ldots, c_s^{n_p} \right)\in \mathbb{R}^{n_p}$.
Unlike SORA, the reliability indices in \autoref{eq:tsm_mc2} are now computed simultaneously and hence involves no approximation. 

\noindent \textbf{Remark 1: }
It is to be noted that extending SORA, similar to Eqs. (\ref{eq:tsm_mc1}) and (\ref{eq:tsm_mc2}), will be difficult as each shift-vector can be high-dimensional. 
On the contrary, the shifted-thresholds of TSM is only one-dimensional and hence it is possible to compute the shifted-thresholds simultaneously.

\noindent \textbf{Remark 2: }
The design-vector in \autoref{eq:tsm_mc2} can become high dimensional if the problem under consideration has several probabilistic-constraints.
This is one of the primary disadvantages of the proposed approach.

To address the limitation of TSM described in remark 2, we propose to carry out the analysis based on the active constraints.
Assuming $n_{ac}$ to be the number of active constraints, \autoref{eq:tsm_mc2} can be rewritten as
\begin{equation}\label{eq:tsm_mc3}
    \begin{split}
        & \bm c_s^* = \argmin_{\bm c_s} \sum_{i=1}^{n_{ac}}\left[ \beta_{\text{Tar},i} - \beta _ {i}\left(\bm c_s \right) \right]^2, \\
        \text{s.t. }\;\;\;\; & \beta_{\text{Tar},i} - \beta _ {i}\left(\bm c_s \right) \le 0, i=1,\ldots, n_{ac}, \\
    \end{split}
\end{equation}
where $n_p$ is substituted for $n_{ac}$.

Despite the advantages of the TSM highlighted above, implementing TSM in practice is a challenge.
Firstly, although reliability indices can be represented as a function of the shifted-thresholds, the mapping from shifted-threshold, $\bm c_s$ to reliability indices, $\beta$ is unknown.
Secondly, although we argued that working with active constraints will reduce the computational effort, how to compute the active constraints is yet to be discussed.

In the next section, we present a dual-surrogate based efficient approach for addressing all the questions raised above.
\subsection{Dual-surrogate based TSM}
Having described the fundamentals of TSM, in this section we present a practical approach for implementing TSM.
The proposed method utilizes two surrogate models and hence, we refer to it as the dual-surrogate based TSM (DS-TSM).
To describe DS-TSM, we revisit \autoref{eq:tsm_mc1}.
Note that even the deterministic optimization involves a constraint of the form $g_i\left(\bm d, \bm {\mu_X}\right) \le c_s^i$.
For a structural engineering problem, the constraints may be of stress or displacement-type and hence, a finite element (FE) solver is required.
If the problem under consideration is complex, a single run of the FE solver takes significant time. 
Therefore, several calls to the FE solver from an optimization subroutine makes the procedure computationally expensive.
To circumvent this, we propose to utilize a surrogate to approximate the constraints, $g_i$ in \autoref{eq:tsm_mc1}. 
We refer to this surrogate as {\it surrogate 1} and denote it as  $\mathcal{M}_1$.
Introducing {\it surrogate 1}, the optimization problem in \autoref{eq:tsm_mc1} can be written as:
\begin{equation}\label{eq:ds-tsm_mc1}
    \begin{split}
        & \bm d^* = \argmin c\left(\bm d \right), \\
        \text{s.t. }\;\;\;\; & \hat g_i\left(\bm d, \bm {\mu_X} \right) \le c_s^i, i=1,\ldots, n_p,\\
        & \bm d^L \le \bm d \le  \bm d^U,
    \end{split}
\end{equation}
where $\hat g_i$ represents an approximation to the actual function $g_i$, obtained by using a surrogate.
This enables us to solve the optimization problem of TSM in an efficient manner.

The next challenge in TSM is to solve \autoref{eq:tsm_mc3}.
Here, we face two challenges - (a) how to learn the mapping from $\bm c_s$ to $\beta$ and (b) how to learn the active constraints.
We tackle the challenges one at a time.
First, we train a surrogate model to learn the mapping from $\bm c_s$ to $\beta$. 
We address this surrogate as {\it surrogate 2} and denote it as $\mathcal{M}_2$.
For training {\it surrogate 2}, we first generate samples for training the surrogate model and then, train a surrogate model to map from $\bm c_s$ to $\beta$. Once trained, the surrogate model can map $\bm c_s$ to $\beta$ and solve the optimization problem in \autoref{eq:tsm_mc3}, where we replace $\beta_i\left( \bm c_s \right)$ with $\hat{\beta}_i\left( \bm c_s \right)$ (surrogate predicted). 

For generating training samples for {\it surrogate 2}, it is essential to first fix the bounds of $c_s^i, \forall i$. 
Assuming that the constraint $g_i$ in \autoref{eq:tsm_mc1} are all of {\it less than type}, it is obvious that $ \left(c_s^i\right)^* \le c_{s,0}^i, \forall i$, where $c_{s,0}^i$ represents the threshold corresponding to the actual probabilistic constraint {\it i.e.}
\begin{align}
    \mathbb{P}\left(G_i \left( \bm d, \bm X \right) > c_{s,0}^i \right) \le P_f^{\text{Tar}},
\end{align}
where 
$P_f^{\text{Tar}}$ is the target-probability of failure.
Therefore, the upper bound of $c_s^i$ is set to $c_{s,0}^i$.
However, no such relation for the lower-bound exists.
In this study, an empirical lower-bound of the form $\alpha c_{s,o}^i$, where $0 < \alpha < 1$, is considered.
Therefore,
\begin{equation}\label{eq:lb}
    \alpha c_{s,0}^i \le c_s^i \le c_{s,0}^i\;\; \; \forall i.
\end{equation}
As a special case, if $c_{s,0}^i = 0$, we set the lower-bound to be $-3$.
Algorithm \ref{alg:bounds} shows the steps involved in setting the bounds for the threshold.

\begin{algorithm}
\caption{Determining the bounds of the thresholds}
\label{alg:bounds}
\textbf{Initialize: }
Provide $c_{s,0}^i$ and $\alpha$.\\
${c}_{U,s}^i \leftarrow c_{s,0}^i$. \\ 
\eIf{$c_{s,0}^i > 0$}{${c}_{L,s}^i \leftarrow \alpha \times c_{s,0}^i$}{
${c}_{L,s}^i \leftarrow -3.0$}
\textbf{Output: }${c}_{U,s}^i, {c}_{L,s}^i$.
\end{algorithm}

\noindent \textbf{Remark 3: }
It is possible that the optimum threshold $\left(c_s^i\right)^*$ may not reside within the bounds obtained from \autoref{alg:bounds}.
However, this will not be an issue as the bounds are only utilized for training the surrogate.
Moreover, the surrogate will also be updated continuously (discussed later) and hence, the accuracy of surrogate will also not be an issue.

Once the input training samples are generated, within the bounds, using \autoref{alg:bounds}, we compute $\beta$ corresponding to the training samples.
To that end, we first solve the optimization problem in \autoref{eq:ds-tsm_mc1} to obtain the design vector $\bm d$. 
Subsequently, we solve a reliability analysis problem to $\beta$.
The steps involved are shown in \autoref{alg:training_beta_generation}.
\begin{algorithm}
\caption{Generating training samples for {\it surrogate 2}}
\label{alg:training_beta_generation}
\textbf{Pre-requisite}:
Trained surrogate 1 ($\mathcal{M}_1$), input training samples $\left(\mathbf{C_s} = \left[\bm c_s^{\left(1\right)}, \ldots, \bm c_s^{\left(N_{s2}\right)} \right]^T\right)$.\\
\For{$i=1,\ldots,N_{s2}$}{
$\bm c_s = \bm c_s^{\left(i\right)}$. \\
Solve the optimization problem in \autoref{eq:ds-tsm_mc1} to compute the corresponding optimal design variables $\bm d$. \\
\For{$j= 1, \ldots, n_p$}{
Using $\mathcal{M}_1$, solve a reliability analysis to compute $\mathbb{P}\left( G_j\left( \bm d, \bm X \right)\right) \le c_{s,0}^j$ and corresponding reliability indices $\beta_j$.\\
}
$\bm{\beta}^{\left(i\right)} = \left[\beta_1, \ldots, \beta_{n_p} \right]$
}
$\mathbf{B}= \left[\bm{\beta}^{\left(1\right)}, \ldots, \bm{\beta}^{\left(N_{s2}\right)} \right]^T$.\\
\textbf{Output:}
$\mathbf{B}$
\end{algorithm}
Once the training samples $\mathbf{C}_s$ and $\mathbf{B}$ have been generated, we can easily train $\mathcal{M}_2$.

We now shift our focus on the other issue {\it i.e.,} determination of active constraints.
In this work, we propose a simple strategy for determining the active constraints.
The strategy is to track the variability of reliability index $\beta$ for each of the constraints to decide whether that constraint is active of not.
The $i$-th constraint is considered to be inactive if $\min\left(\mathbf{B}_{:,i}\right) \ge {\epsilon_1}\beta_{\text{Tar}}^i$, where
$\beta_{\text{Tar}}^i$ is the target reliability index of the $i$-th constraint and $\mathbf{B}_{:,i}$ represents the
$i$-th column of matrix $\mathbf{B}$ in \autoref{alg:training_beta_generation} and $\epsilon_1$ is a positive factor considered to be 2.5.
Algorithm \ref{alg:active_cons} depicts the steps involved in selecting the active constraints.

\begin{algorithm}
\caption{Selection of active constraints}
\label{alg:active_cons}
\textbf{Pre-requisite: }$\bm{\beta}_{\text{Tar}}$, $\mathbf{B}$, $\epsilon_1$.\\
\textbf{Set: }Initiate a null vector \texttt{index}\\
\For {$i=1,\ldots, n_p$}{
$\beta_{min} \leftarrow \min\left( \mathbf{B}_{:,i}\right)$.\\
\If{$\beta_{min} \leq \epsilon_1 \times {\beta}_{Tar}^i$}{
Augment $i$ to the vector \texttt{index}.
}
}
\textbf{Output: } \texttt{index}.
\end{algorithm}

The primary advantage of this active-constraint selection technique resides in the fact 
that it can be viewed as a simple extension to \autoref{alg:training_beta_generation}.
In other words, no additional computational overhead is associated with this method.

Finally, we focus on the DS-TSM.
In this algorithm we utilize algorithms \ref{alg:bounds} -- \ref{alg:active_cons} to iterate over Eqs. (\ref{eq:tsm_mc3}) and (\ref{eq:ds-tsm_mc1}).
At each iteration, $\mathcal{M}_2$ is updated based on results obtained from the immediate previous step.
The algorithm stops when convergence criteria, defined in terms of $\beta$ is satisfied.
The steps involved in the DS-TSM are shown in \autoref{alg:threshold_shift_method}.

\begin{algorithm}
\caption{Dual-surrogate based threshold shift method}
\label{alg:threshold_shift_method}
\textbf{Input: }
$\bm d_l$, $\bm d_u$, $\bm c_{s,0}$, $\bm{\beta}_{\text{Tar}}$, $\alpha$, $\epsilon_1$, $\epsilon_{tol}$. \\
Generate training samples for $\mathcal{M}_1$. \\
Train $\mathcal{M}_1$.\\
Using $\alpha$, obtain bounds of the constraint thresholds by using \autoref{alg:bounds}.\\
Generate $\mathbf{C}_s$ based on the generated bounds. \\
Generate $\mathbf{B}$ by using \autoref{alg:training_beta_generation}.\\
Based on $\epsilon_1$, identify the active constraints, \texttt{index} by using \autoref{alg:active_cons}. \\
$\mathbf{C}_{s,red} \leftarrow \mathbf{C}_s\left(\right.$\texttt{index}$\left.\right)$, $\mathbf{B}_{red} \leftarrow \mathbf{B}_{red}\left(\right.$\texttt{index}$\left.\right)$.\\
\textbf{Set } $\bm{c}_s \leftarrow \bm c_{s,0}$. \\
\Repeat{$ \min \left| \bm{\beta}_{\text{Tar}} - \bm{\beta}^*  \right| \le \epsilon_{tol}$ and $\max \left( \bm{\beta}_{\text{Tar}} - \bm{\beta}^*  \right) \le 0$}{
Train $\mathcal{M}_2$ to using $\mathbf{C}_{s,red}$ and $\mathbf{B}_{red}$.\\
Obtain $\bm c_s^*$ by solving \autoref{eq:tsm_mc3}.\\
$\bm{c}_s\left(\right.$\texttt{index}$\left.\right) \leftarrow \bm c_s^*$.\\
With $\bm c_s^*$ as threshold shift, $\bm d_l$ as lower-bound and $\bm d_u$ as upper bound, obtain $\bm d$ by solving \autoref{eq:ds-tsm_mc1}. \\
Perform reliability analysis corresponding to the active constraints to obtain $\bm{\beta}^*$.\\
$\mathbf{B}_{red} = \left[\mathbf{B}_{red}; \bm{\beta}^* \right]$, $\mathbf{C}_{s,red} = \left[\mathbf{C}_{s,red}; \bm{c}_s^* \right]$.\\
}
$\bm d^* = \bm d$. \\
\textbf{Output: }
$\bm d^*$.
\end{algorithm}

\begin{figure*}[htbp!]
    \centering
    \includegraphics[width=\textwidth]{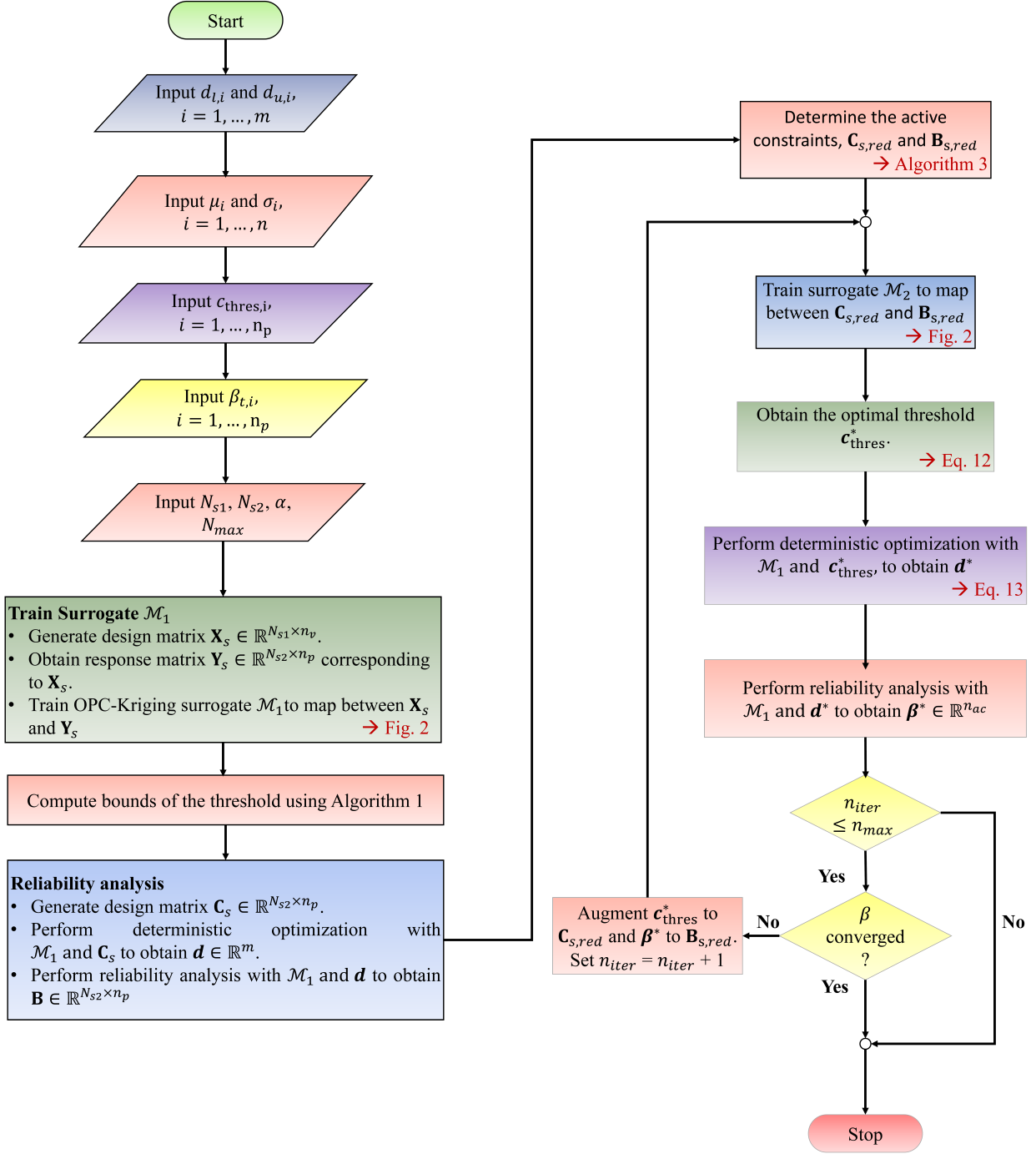}
    \caption{Optimal PC-Kriging based threshold shift method}
    \label{fig:opcK_tsm}
\end{figure*}

From \autoref{alg:threshold_shift_method}, it is obvious that the performance of the proposed DS-TSM will depend on the performance of the two surrogates, $\mathcal{M}_1$ and $\mathcal{M}_2$.
Therefore, it is essential to use an accurate surrogate modeling technique and tune it in the best possible manner.
To that end, an efficient surrogate, referred to as polynomial chaos based Kriging (PC-Kriging) has been used.
Both $\mathcal{M}_1$ and $\mathcal{M}_2$ are trained PC-Kriging models. 
A brief description on PC-Kriging model is provided in Appendix \ref{appendix-pc_kriging}.
The overall framework of OPC-Kriging based TSM is shown in \autoref{fig:opcK_tsm}.

\noindent \textbf{Remark 4: }
Note that the second surrogate ({\it surrogate 2}) within the proposed DS-TSM is trained based on training samples generated using the first surrogate (optimization and reliability analysis performed with the help of first surrogate).
Therefore, the computational cost associated with the proposed approach is equivalent to the number of training samples required for training the first surrogate. 
\section{Numerical examples}
\label{sec:numericals}
In this section, six numerical examples are presented to illustrate the performance of the proposed approach. 
The examples are well-established benchmark problems available from literature.
While the first two problems are mathematical in nature, the last four problems are physical problems involving engineering systems. 
For all the problems, comparison with results available in literature is presented.

In all the problems presented in this section, OPC-Kriging has been formulated by considering maximum degree to be 5. 
For the Kriging part in OPC-Kriging, {\it mattern} correlation function has been used due to its already proven superiority over most of the other covariance functions.
Unless otherwise mentioned, Sobol sequence \cite{Ref_paper3_47} has been used for generating the training samples.
Number of training samples varies in each problem and is mentioned while defining the problem.
We have considered $\alpha = \frac{1}{3}$, $\epsilon_1 = 2.5$ and $\epsilon_{\text{tol}} = 0.001$. Maximum number of iterations allowed is 100.

The proposed TSM involves three major modifications -- (a) computing all the thresholds by solving a single optimization algorithm, (b) using OPC-Kriging to replace the actual limit state function and (c) operating only based on the active probabilistic constraints. 
Hereafter, we refer to the three modfications discussed above as {\it Modification (a), Modification (b)} and {\it Modification (c)}, respectively.
Results illustrating the advantages (both in terms of accuracy and efficiency) gained by using each of these modifications have also been discussed.

\subsection{Two-dimensional non-linear limit state surface}
As the first example, we consider a RBDO problem with two-dimensional non-linear limit state surface.
This example has been previously studied in \cite{Lee2008}. 
The limit state functions are:
\begin{subequations}\label{seq:ls_prob1}
\begin{equation}
    c_1\left(\bm x \right) = -x_1 \sin \left(4x_1\right) - 1.1x_2\sin \left(2x_2\right),
\end{equation}
\begin{equation}
    c_2\left(\bm x \right) = x_1 + x_2 - 3,
\end{equation}
\end{subequations}
where $\bm x = \left[x_1, x_2 \right]$ are the two independent normal random variables associated with this problem.
We consider $\sigma_{x_1} = \sigma_{x_2} = 0.1$, where $\sigma_i$ represents standard deviation (std. dev.) of variable $i$.
The RBDO problem reads as:
\begin{equation}\label{eq:prob1_ps}
    \begin{split}
        \text{Minimize } & c\left( \bm d \right) = \left(d_1 - 3.7 \right)^2 + \left( d_2 - 4 \right)^2, \\
        \text{s.t.}\;\;\;\; &\mathbb{P}\left( c_l\left( \bm x \right) \le 0\right) \le \Phi \left( \beta_{\text{Tar}}^l\right), \; l=1,2, \\
        & 0 \le d_1 \le 3.7, \\
        & 0 \le d_2 \le 4.0. \\
    \end{split}
\end{equation}
In \autoref{eq:prob1_ps}, $d_i = \mu_{x_i},\; i=1,2$, where $\mu_i$ represents mean of variable $i$.
$\Phi$ in \autoref{eq:prob1_ps} represents cumulative distribution function of the standard normal variables.
For both the constraints, we consider $\beta_{\text{Tar}}$ to be 2. {\it Surrogate 1} is trained with 32 training samples for $g_1$ and $10$ samples for $g_2$.

\autoref{tab:prob1_res} shows the results obtained using different methods.
For the sake of comparison, results obtained using deterministic design optimization (DDO) have also been presented.
As expected, results obtained using RBDO techniques are conservative as compared to the DDO results.
This is expected as the RBDO methods provide additional insurance against failure due to the violation of the probabilistic constraints.

As for the performance of the RBDO techniques, it is observed that TSM yields the best results followed by Meta-RBDO proposed in \cite{Dubourg2011}.
PMA and PMA w/Kriging \cite{Lee2008} are found to violate the probabilistic constraints and hence, are unsafe.
Also from the efficiency point of view, TSM is found to be the most efficient with 42 actual function evaluations.
Meta-RBDO on the other hand, requires 90 function evaluations.

Finally, we compare the advantages gained by each of the modifications (\autoref{tab:p1_diffmod}). It is observed that no additional advantage is obtained by {\it Modification (a)} (computing each thresholds by solving a separate optimization problem). 
This is because the number of active constraints for this problem is only 1.
The fact that we are using a surrogate model to approximate the limit state function (i.e., {\it Modification b}) provides significant computational efficiency (with almost no compromise on the accuracy).
Finally, without active constraints, the second surrogate model ({\it surrogate 2}) requires more samples to converge.
As a consequence, number of iterations to achieve convergence increases significantly.

\begin{table}[htbp!]
    \centering
    \caption{Results for Problem 1. Results for PMA and PMA w/Kriging are taken from \cite{Lee2008}. Results for for Brute force and Meta-RBDO are taken from \cite{Dubourg2011}. Results obtained using Brute force method are the benchmark solution.}
    \label{tab:prob1_res}
    \begin{tabular}{lccccccc}
    \hline
        \textbf{Methods} & $d_1$ & $d_2$ & $c\left(\bm d \right)$ & $\beta_1$ & $\beta_2$ & Function calls & Opt. iter$^*$  \\ \hline \hline
        \textbf{DDO} & 2.97 & 3.41 & 0.89 & $-0.10$ & $> 8.00$ & 70/70 & 10 \\ \hline
        Brute force & 2.84 & 3.23 & 1.33 & $2.00$  & $> 8.00$ & $\approx 10^7/10^7$ & 10 \\
        PMA & 2.82 & 3.30 & 1.26 & 1.67 & $> 8.00$ & 296 & 7 \\
        PMA w/Kriging & 2.82 & 3.30 & 1.26 & 1.67 & $> 8.00$ & 90 & 7 \\
        Meta-RBDO & 2.81 & 3.25 & 1.35 & 2.00 & $> 8.00$ & 90 & 10 \\
        DS-TSM & 2.86 & 3.21 & 1.33 & 2.00 & $>8.00$ & 42 & 10 \\ \hline
        \multicolumn{4}{l}{$^*$ Number of iterations for optimization.}
    \end{tabular}
\end{table}

\begin{table}[htbp!]
    \centering
    \caption{Advantage gained using the three modifications of TSM.}
    \label{tab:p1_diffmod}
    \begin{tabular}{lccc}
    \hline
    \textbf{Methods} & $c\left(\bm d \right)$ & Function calls & Opt. iter \\ \hline \hline 
    \textbf{TSM} & 1.33 & 42 & 10 \\ \hline
    TSM (without (a))$^*$ & 1.33 & 42 & 10 \\ 
    TSM (without (b))$^\#$ & 1.33 & $10^6/10^6$ & 10 \\
    TSM (without (c))$^{\dagger}$ & 1.33 & 42 & 25 \\ \hline
    \multicolumn{4}{l}{$^*$Thresholds for each constraint computed separately.}\\
    \multicolumn{4}{l}{$^\#$No surrogate model is used to approximate the limit state function.} \\
    \multicolumn{4}{l}{$^{\dagger}$Active constraint algorithm has not been utilized.}
    \end{tabular}
\end{table}
\subsection{Highly non-linear problem}\label{subsec:highly_non-linear}
In this example a highly non-linear problem has been considered. This problem has previously been studied in \cite{Hao2017} and \cite{Keshtegar2018}. The problem statement reads as:
\begin{equation}\label{eq:prob3_ps}
    \begin{split}
        \text{Minimize } & f\left( \bm d \right) = -\frac{(d_1+d_2-10)^2}{30} - \frac{(d_1-d_2+10)^2}{120}, \\
        \text{s.t.}\;\;\;\; &\mathbb{P}\left(g_j(\left( \bm x \right) > 0\right) \le \Phi \left( -\beta_{\text{Tar}}^{j},\right) \; j=1,2,3, \\
    \end{split}
\end{equation}
where
\begin{equation}\label{eq:prob3_constraints}
    \begin{split}
        g_1 &= 1-\frac{x_1^2x_2}{20}\\
        g_2 &= -1+(Y-6)^2 +(Y-6)^3 - 0.6(Y-6)^4+ Z \\
        g_3 &= 1- \frac{80}{x_1^2 + 8x_2+5}\\
        Y &= 0.9063x_1 + 0.4226x_2\\
        Z &= 0.4226x_1 -0.9063x_2\\
        &0 \le d_i \le 10, \; \; \; x_i \approx N\left(d_i, 0.3^2\right) for \;\;i =1,2 \\
        &\bm{d}^0 = [5,5], \beta_{\text{Tar}}^{1} = \beta_{\text{Tar}}^{2} = \beta_{\text{Tar}}^{3} = 3.0,\\
\end{split}
\end{equation}
In \autoref{eq:prob3_constraints}, $x_1$ and $x_2$ are the normal random variables. $d_1$ and $d_2$ are the design variables in \autoref{eq:prob3_ps}. For implementing the proposed approach, {\it Surrogate 1} is trained with 32 training samples for $g_1$, 128 samples for $g_2$ and $64$ samples for $g_3$. Results corresponding to the chaos control methods are obtained from \cite{Keshtegar2018}.
 
\autoref{tab:prob3_res} shows the results obtained using different methods.
For the sake of comparison, results obtained using DDO have also been presented.
As expected, results obtained using RBDO techniques are conservative as compared to the DDO results.

As for the performance of the RBDO techniques, it is observed that the results obtained using TSM are closer to the benchmark solution obtained using brute-fore method. 
All the chaos control methods are found to violate the probabilistic constraints.
From the efficiency point of view also, TSM is found to be the most efficient with 224 actual function evaluations.

Finally, we compare the advantages gained by each of the modifications (\autoref{tab:p3_diffmod}). 
Without {\it Modification (a)}, the results are found to be less accurate.
The number of iterations required is also found to be significantly more.
The difference in the result is due to the fact that for this problem, we have two active constraints.
The fact that we are using a surrogate model to approximate the limit state function provides significant computational efficiency (with slight compromise on the accuracy).
Without this, the number of function evaluations increases significantly.
Finally, without the active constraint determination scheme, more iterations are required to achieve convergence.

\begin{table}[htbp!]
    \centering
    \caption{Results for problem 2. Results corresponding to hybrid chaos control (HCC), adaptive chaos control (ACC), modified mean value approach (MMV) are obtained from \cite{Keshtegar2018}. Results obtained using Brute force method are the benchmark solution.}
    \label{tab:prob3_res}
    \begin{tabular}{ccccccc}
    \hline
    \textbf{Methods} & $d_1$ & $d_2$ & $f\left(\bm d \right)$ & $\beta_{\text{min}}$ & Function calls & Opt. iter. \\ \hline \hline
    \textbf{DDO} & 5.197 & 0.741 & -2.29 & $-0.022$ &  110 & 35 \\ \hline
    Brute Force & 4.689 & 1.908 & -1.747 & 3.0 &  $3 \times 10^7$   & 8 \\
    HCC & 4.558 & 1.964 & -1.724 & 2.95 &  1629 & 10 \\
    ACC & 4.558 & 1.964 & -1.724 & 2.95 &  861 & 10 \\
    MMV & 4.558 & 1.964 & -1.724 & 2.95 &  1593 & 10\\
    TSM & 4.607 & 1.949 & -1.731 & 3.00 & 224 (32+128+64) & 12 \\ \hline 
    \end{tabular}
\end{table}

\begin{table}[htbp!]
    \centering
    \caption{Advantage gained using the three modification of TSM (Problem 2).}
    \label{tab:p3_diffmod}
    \begin{tabular}{lccc}
    \hline
    \textbf{Methods} & $f\left(\bm d \right)$ & Function calls & Opt. iter \\ \hline \hline 
    \textbf{TSM} & -1.731 & 224 & 12 \\ \hline
    TSM (without (a))$^*$ & -1.71 & 224 & 53 \\ 
    TSM (without (b))$^\#$ & -1.747 & $21 \times 10^5$ & 12 \\
    TSM (without (c))$^{\dagger}$ & -1.731 & 224 & 32 \\ \hline
    \multicolumn{4}{l}{$^*$Thresholds for each constraint computed separately.}\\
    \multicolumn{4}{l}{$^\#$No surrogate model is used to approximate the limit state function.} \\
    \multicolumn{4}{l}{$^{\dagger}$Active constraint algorithm has not been utilized.}
    \end{tabular}
\end{table}

\subsection{Automobile with front I-beam axle}\label{subsec:I-Beam_axle}
\begin{figure}[htbp!]
    \centering
    \includegraphics[width = 0.5\textwidth]{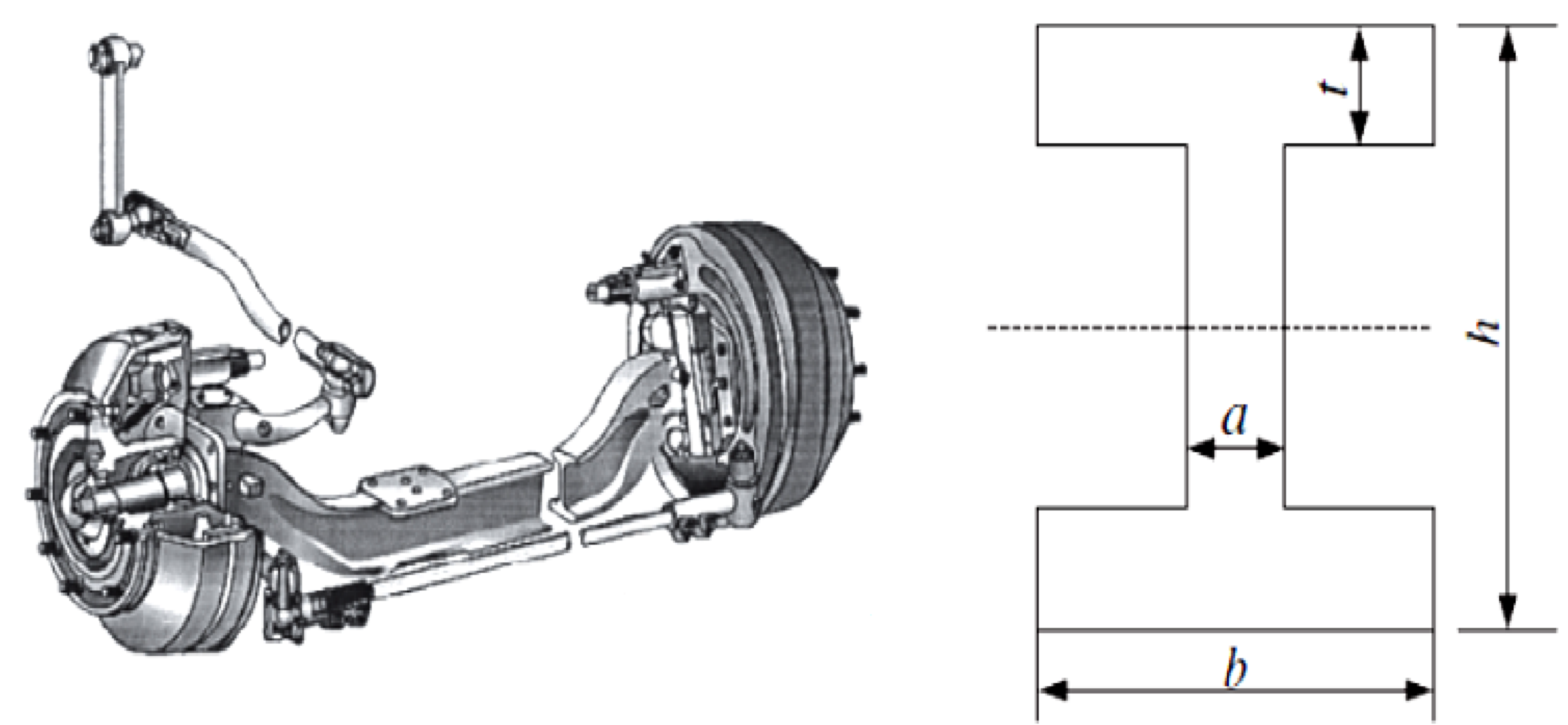}
    \caption{Geometrical description and loading conditions for the automobile with front I-beam axle.}
    \label{fig:prob4}
\end{figure}

In this problem we consider the front I-beam axle of an automobile. \autoref{fig:prob4} shows a diagrammatic representation of the problem. This problem has previously been studied in \cite{Keshtegar2018SC, Keshtegar2018}. The optimization problem reads as: 
\begin{equation}\label{eq:prob4_ps}
    \begin{split}
        \text{Minimize } & f\left( \bm d \right) = ha + 2(b-a)t, \\
        \text{s.t.}\;\;\;\; &\mathbb{P}\left(g\left( \bm x \right) > 0\right) \le \Phi \left( -\beta_{\text{Tar}}\right) \\
        \text{where}\;\;\;\; & g = \sqrt{\sigma _m ^2 +3 \tau ^2} -\sigma _s\\
        & \sigma _m = \frac{M}{W_x}, \; \; \tau = \frac{T}{W_{\rho}} \\
        & W_x = \frac{a(h-2t)^3}{6h} + \frac{b}{6h}[h^3 - (h-2t)^3]\\
        & W_{\rho} = 0.8bt^2 + 0.4\frac{a^3(h-2t)}{t}\\
        & 60 \le h \le 100, \; \; \; 10 \le a \le 20 \\
        & 70 \le b \le 120, \; \; \; 10\le t\le 20 \\
        & [a,b,t,h]^0 = [85, 12,75,12]\\ 
        & \beta_{\text{Tar}} = 3.0, \\
    \end{split}
\end{equation}
where $\sigma _s = 460$ MPa is the yielding stress of the material. $\tau$ and $\sigma_m$ are the maximum shear stress and the normal stress developed due to torque ($\tau$) and bending moment ($M$), respectively. $W_x$ and $W_{\rho}$ are section factor and polar section factor \cite{Keshtegar2018SC}, respectively. This example has four design random variables; $a$, $b$, $t$ and $h$. The statistical parameters of the random variables are provided in \autoref{tab:I-beam_axle}.
For utilizing the proposed approach to solve this problem, {\it Surrogate 1} has been trained with 128 training samples.
Results corresponding to the chaos control methods are obtained from \cite{Keshtegar2018}.
\begin{table}[htbp!]
    \centering
    \caption{Statistical parameters of random parameters in \autoref{eq:prob4_ps}.}
    \label{tab:I-beam_axle}
    \begin{tabular}{ccccc}
    \hline
    \text{Random variables} & \text{Unit} & \text{Type} & \text{Mean} & \text{Std.dev.}\\ 
    \hline
    $a$ & mm & Normal & - &  0.060 \\
    $b$ & mm & Normal & - &  0.325 \\
    $t$ & mm & Normal & - &  0.070 \\
    $h$ & mm & Normal & - &  0.425 \\
    $M$ & N-mm & Normal & $3.5\times10^6$ &  $1.75\times10^5$ \\
    $T$ & N-mm & Normal & $3.1\times10^6$ &  $1.55\times10^5$ \\
    \hline
    \end{tabular}
\end{table}

\autoref{tab:prob4_res} shows the results obtained using different methods.
It is observed that TSM yields the best results (lower objective function and no constraint violation).
From the efficiency point of view also, TSM is found to be the most efficient with 128 actual function evaluations.

Finally, we compare the advantages gained by each of the modifications (\autoref{tab:p4_diffmod}). 
Since, we have only one probabilistic constraint, {\it Modifications (a) and (c)} have no effect.
On the other hand, {\it Modification (b)} provides significant computational efficiency, and without this modification, the number of function evaluations increases significantly.

\begin{table}[htbp!]
    \centering
    \caption{Results for problem 3. Results corresponding to HCC, ACC, step-modified chaos control (SMCC) are obtained from \cite{Keshtegar2018}. Results obtained using Brute force method are the benchmark solution.}
    \label{tab:prob4_res}
    \begin{tabular}{ccccc}
    \hline
    \textbf{Methods} & HCC & ACC & SMCC & TSM \\ \hline
    $t^*$ & 16.497 & 16.497 & 16.497 & 16.303 \\
    $h^*$ & 60 & 60 & 60 & 60 \\
    $a^*$ & 10 & 10 & 10 & 10 \\
    $b^*$ & 70 & 70 & 70 & 70 \\
    $f\left(\bm d \right)$ & 2024.832 & 2024.832 & 2024.832 & 2015.144 \\
    Opt. iter. & 6 & 6 & 6 & 8 \\
    func. calls & 1303 & 415 & 1637 & 128 \\ \hline
    \end{tabular}
\end{table}

\begin{table}[htbp!]
    \centering
    \caption{Advantage gained using the three modification of TSM (Problem 3).}
    \label{tab:p4_diffmod}
    \begin{tabular}{lccc}
    \hline
    \textbf{Methods} & $f\left(\bm d \right)$ & Function calls & Opt. iter \\ \hline \hline
    \textbf{TSM} & 2015.144 & 128 & 8 \\ \hline
    TSM (without (a))$^*$ & 2015.144 & 128 & 8 \\ 
    TSM (without (b))$^\#$ & 2015.144 & $8 \times 10^5$ & 8 \\
    TSM (without (c))$^{\dagger}$ & 2015.144 & 128 & 8 \\ \hline
    \multicolumn{4}{l}{$^*$thresholds for each constraint computed separately}\\
    \multicolumn{4}{l}{$^\#$No surrogate model is used to approximate the limit state function} \\
    \multicolumn{4}{l}{$^{\dagger}$Active constraint algorithm has not been utilized}
    \end{tabular}
\end{table}

\subsection{Bracket structure}\label{subsec:bracket_str}
\begin{figure}[htbp!]
    \centering
    \includegraphics[width = 0.4\textwidth]{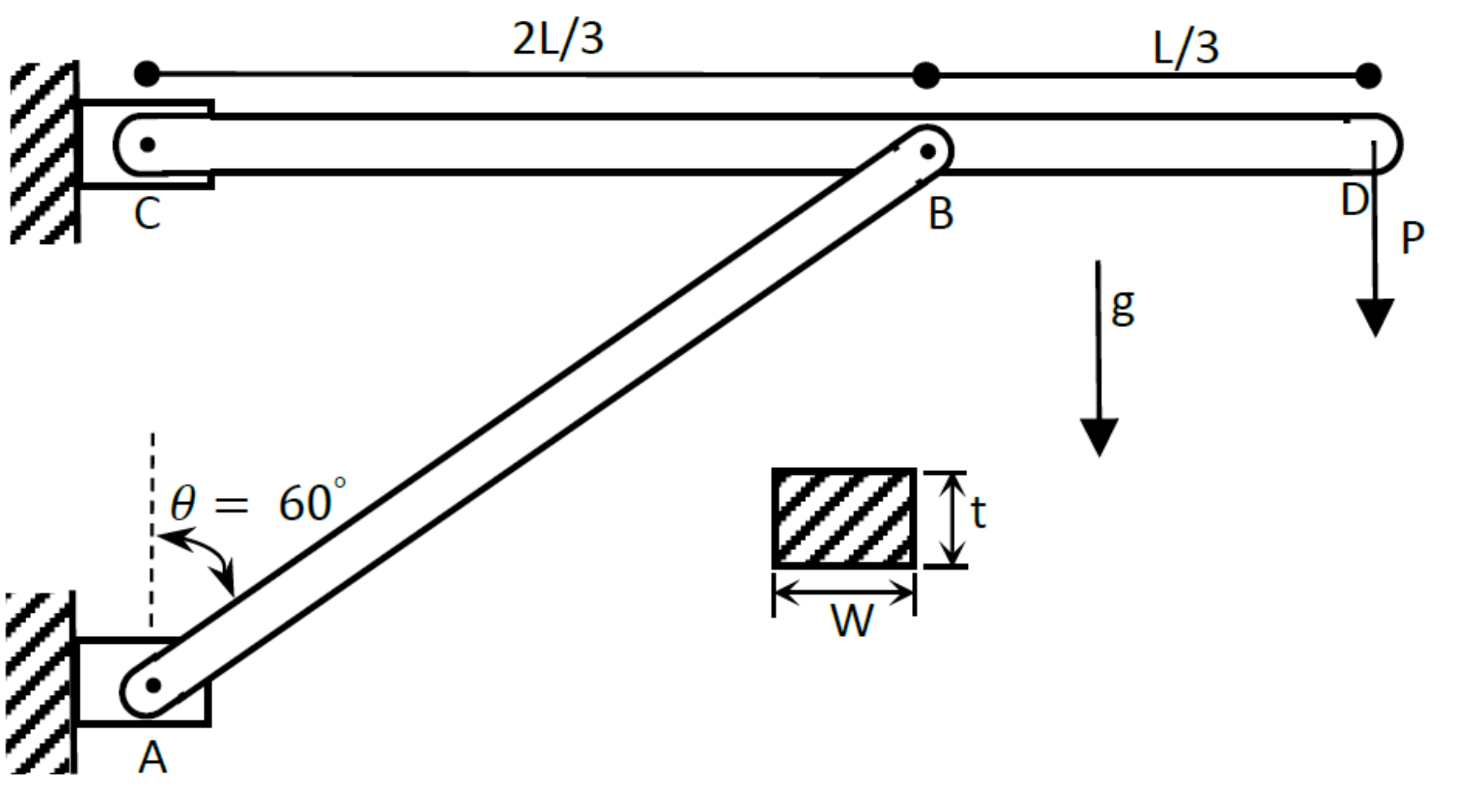}
    \caption{Geometrical description and loading conditions for the structural bracket.}
    \label{fig:prob6}
\end{figure}

In this example, we consider a bracket structure as shown in \autoref{fig:prob6}. It is loaded by its own self weight and by a vertical load. Both members have the same thickness, $t$ but have different widths, $W_{AB}$ and $W_{CD}$.
The RBDO problem seeks to find the rectangular cross-sections of the two structural members that minimize the expected overall structural weight, which is computed as
follows:
\begin{equation}\label{eq:prob6_ps}
    \begin{split}
        \text{Minimize } & f\left( \bm d \right) = \rho tL\left(\frac{4\sqrt(3)}{9}W_{AB} + W_{CD}\right), \\
        &\mathbb{P}\left(g_j(\left( \bm x \right) < 0\right) \le \Phi \left( -\beta_{\text{Tar}}^{j},\right) \; j=1,2, \\
        \text{where } & g_1 = \sigma_y - \sigma_B(\bm{x})\\
         & g_2 = F_{buckling}(\bm{x}) - F_{AB}(\bm{x})\\
         & \sigma_B(\bm{x}) = \frac{6M_{B}}{W_{CD}t^2}\;\;  \\&\text{with}\;\; M_B = \frac{PL}{3} + \frac{\rho gW_{CD}tL^2}{18}\\
         & F_{AB}(\bm{x}) = \frac{1}{\cos \theta}\left(\frac{3P}{2} + \frac{3\rho g W_{CD}tL}{4}\right)\\
         & F_{buckling} = \frac{9\pi ^2EtW_{AB}^3\sin ^2 \theta}{48L^2}\\
         & 50 \le W_{AB}, W_{CD}, t \le 300\;mm\\
         & \beta_{\text{Tar}}^{j} = 2.0\; j=1,2. 
    \end{split}
\end{equation}
The statistical parameters of the random variables are provided in \autoref{tab:bracket_str}, where $COV$ represents the coefficient of variance.
\begin{table}[htbp!]
    \centering
    \caption{Statistical parameters of random parameters in \autoref{eq:prob6_ps}.}
    \label{tab:bracket_str}
    \begin{tabular}{ccccc}
    \hline
    \text{Random variables} & \text{Unit} & \text{Type} & \text{Mean} & \text{COV}.\\ 
    \hline
    $P$ & kN & Gumbel & 100 &  15\% \\
    $E$ & GPa & Gumbel & 200 &  8\% \\
    $\sigma_y$ & MPa & Lognormal & 225 &  8\% \\
    $\rho$ & $kg/m^3$ & Weibull & 7860 &  10\% \\
    $L$ & m & Normal & 5 & 5\% \\
    $W_{AB}$ & mm & Normal & - &  5\% \\
    $W_{CD}$ & mm & Normal & - &  5\% \\
    $t$ & mm & Normal & - &  5\% \\
    \hline
    \end{tabular}
\end{table}
This problem has previously been studied in \cite{Dubourg2011}.
For using TSM to solve this problem, {\it Surrogate 1} has been trained using 128 samples for $g_1$ and 64 samples for $g_2$. 

\autoref{tab:prob4_res} shows the results obtained using different methods.
For the sake of comparison, results obtained using DDO with the partial safety factor has also been presented.
Compared to the RBDO results obtained using various methods, it is observed that DDO with partial safety factor yields highly conservative results.

As for the performance of the RBDO techniques, it is observed that TSM yields the best results followed by meta-RBDO. 
Results corresponding to SORA and nested-FORM are obtained from \cite{Du2004}. 
The meta-RBDO results are obtained from \cite{Dubourg2011}.
As for efficiency, while TSM requires 172 actual function evaluations, meta-RBDO needs 250 function evaluations. 
Both SORA and nested-FORM require 1,675 function evaluations.

Finally, we compare the advantages gained by each of the modifications (\autoref{tab:p5_diffmod}). 
Without {\it Modification (a)}, the results are found to be less accurate.
The number of iterations required is also found to be significantly more.
The difference in the results is due to the fact that for this problem, both the constraints are active (i.e., no inactive constraint).
The fact that we are using a surrogate model to approximate the limit state function provides significant computational efficiency Without the surrogate model ({\it Surrogate 1}), the number of function evaluations increases significantly.
The last modification has no effect for this problem as there exists no inactive constraint.

\begin{table}[htbp!]
    \centering
    \caption{Results for the bracket structure. Results corresponding to Meta-RBDO is obtained from \cite{Dubourg2011}. Results corresponding to DDO, SORA and Nested FORM are obtained from \cite{Du2004}.}
    \label{tab:prob5_res}
    \begin{tabular}{llccc}
    \hline
    \textbf{Method} & \textbf{Opt. design} (mm) & \textbf{Cost} (kg) & \textbf{func. calls} & \textbf{Reliability} \\ \hline
    \multirow{3}{*}{\textbf{DDO w/PSF}} & $w_{AB} = 61$ & \multirow{3}{*}{2,632} & \multirow{3}{*}{40} & \multirow{2}{*}{$\beta^1$ = 4.83} \\
    & $w_{CD} = 202$ & & & \multirow{2}{*}{$\beta^2 = 2.83$} \\
    & $t = 269$ & & & \\ \hline
    \multirow{3}{*}{SORA} & $w_{AB} = 61$ &  \multirow{3}{*}{1,675} & \multirow{3}{*}{1,340} & \multirow{2}{*}{$\beta^1 = 1.96$} \\
    & $w_{CD} = 157$ & & & \multirow{2}{*}{$\beta^2 = 1.98$} \\
    & $t = 209$ & & & \\
    & & & & \\
    \multirow{3}{*}{Nested-FORM} & $w_{AB} = 61$ &  \multirow{3}{*}{1,675} & \multirow{3}{*}{2,340} & \multirow{2}{*}{$\beta^1 = 1.96$} \\
    & $w_{CD} = 157$ & & & \multirow{2}{*}{$\beta^2 = 1.98$} \\
    & $t = 209$ & & & \\
    & & & & \\
    \multirow{3}{*}{Meta-RBDO} & $w_{AB} = 58$ &  \multirow{3}{*}{1,584} & \multirow{2}{*}{160} & \multirow{2}{*}{$\beta^1 = 1.98$} \\
    & $w_{CD} = 128$ &  & \multirow{2}{*}{90} & \multirow{2}{*}{$\beta^2 = 1.94$} \\
    & $t = 209$ & & & \\
    & & & & \\
    \multirow{3}{*}{TSM} & $w_{AB} = 54.75$ &  \multirow{3}{*}{1,495} & \multirow{2}{*}{128} & \multirow{2}{*}{$\beta^1 = 2.41$} \\
    & $w_{CD} = 84.715$ & & \multirow{2}{*}{64} & \multirow{2}{*}{$\beta^2 = 2.19$} \\
    & $t = 300$ & & & \\
    \hline
    \end{tabular}
\end{table}

\begin{table}[htbp!]
    \centering
    \caption{Advantage gained using the three modification of TSM (Problem 4).}
    \label{tab:p5_diffmod}
    \begin{tabular}{lccc}
    \hline
    \textbf{Methods} & $f\left(\bm d \right)$ & Function calls & Opt. iter \\ \hline \hline
    \textbf{TSM} & 1,495 & 172 & 40 \\ \hline
    TSM (without (a))$^*$ & 1,637 & 172 & 100 \\ 
    TSM (without (b))$^\#$ & 1,492 & $3 \times  10^6$ & 32 \\
    TSM (without (c))$^{\dagger}$ & 1,495 & 172 & 40 \\ \hline
    \multicolumn{4}{l}{$^*$thresholds for each constraint computed separately}\\
    \multicolumn{4}{l}{$^\#$No surrogate model is used to approximate the limit state function} \\
    \multicolumn{4}{l}{$^{\dagger}$Active constraint algorithm has not been utilized}
    \end{tabular}
\end{table}

\subsection{Roof truss}\label{subsec:roof_truss}
\begin{figure}[htbp!]
    \centering
    \includegraphics[width = 0.5\textwidth]{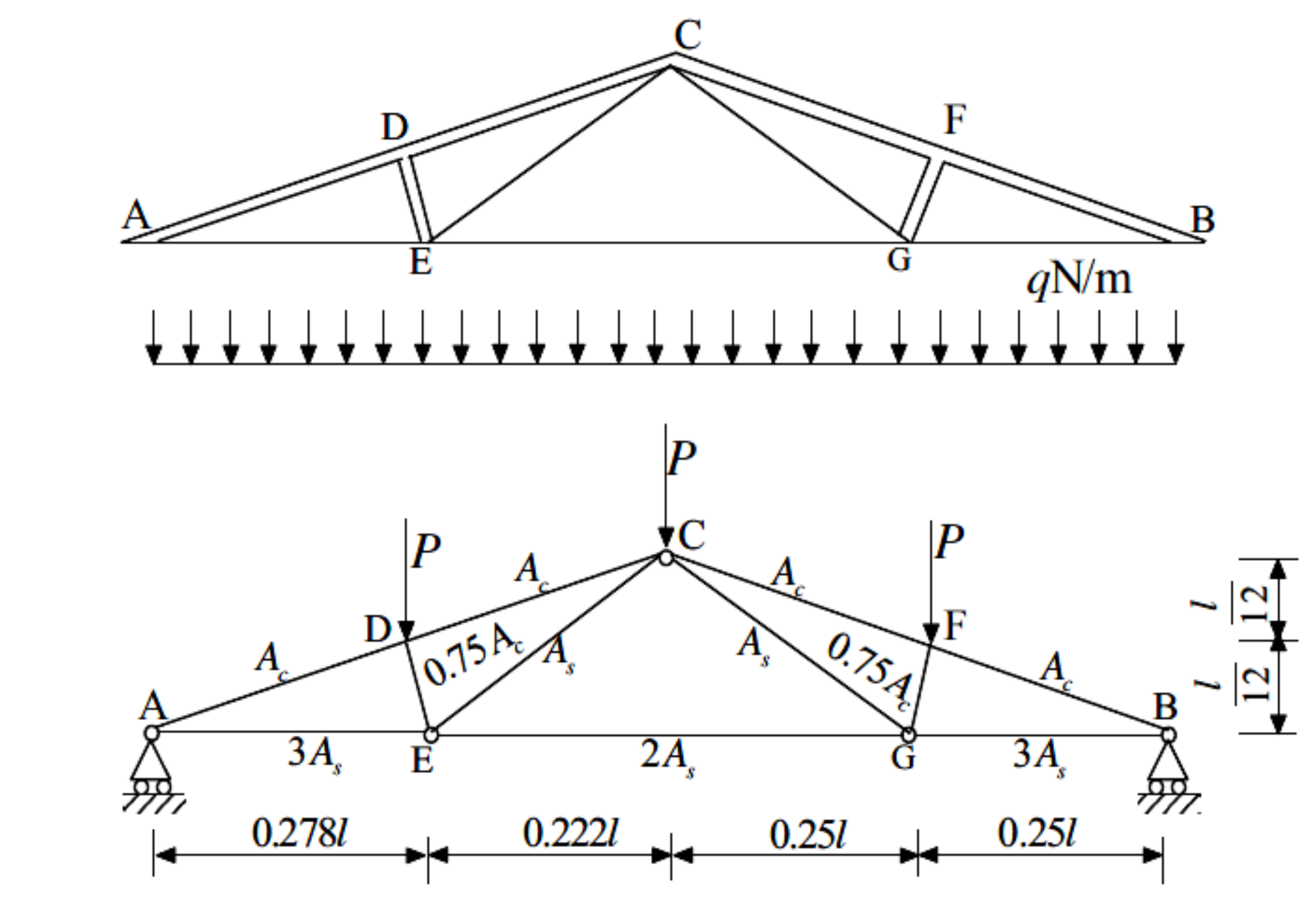}
    \caption{Geometrical description and loading conditions for the roof truss.}
    \label{fig:prob5}
\end{figure}

A roof truss with uniform loading is considered in this example. The problem setup is shown in \autoref{fig:prob5}. The top compression members are made of reinforced concrete while the bottom tensile members are made of steel. The RBDO model based on the reinforced concrete and steel components is given as follows:
\begin{equation}\label{eq:prob5_ps}
    \begin{split}
        \text{Minimize } & f\left( \bm d \right) = 20224 A_s +364 A_C, \\
        \text{s.t.}\;\;\;\; &\mathbb{P}\left(g(\left( \bm x \right) = 0.03 - \frac{ql^2}{2}\left(\frac{3.81}{A_cE_c} + \frac{1.13}{A_s E_s}\right) \le 0\right)\\
        &\le \Phi \left( -\beta_{\text{Tar}}\right) \\
        \text{where}\;\;\;\; &0.0006 \le A_s \le 0.0012,  \;\; 0.018 \le A_c\le 0.063\\
        &[A_s, A_c]^0 = [0.001, 0.042], \beta_{\text{Tar}} = 3.0,\\
    \end{split}
\end{equation}
where $A_s$ and $A_c$ are random design variables and $q$, $l$, $E_s$ and $E_c$ are normal random variables. The statistical properties of the random variables are listed in \autoref{tab:roof_truss}.
For using the proposed approach to solve this problem, we have utilized 164 samples to train {\it Surrogate 1}. Results corresponding to the chaos control methods are obtained from \cite{Keshtegar2018}.
\begin{table}[htbp!]
    \centering
    \caption{Statistical parameters of random parameters in \autoref{eq:prob5_ps}.}
    \label{tab:roof_truss}
    \begin{tabular}{ccccc}
    \hline
    \text{Random variables} & \text{units} &\text{Mean} & \text{Std.dev.} & \text{Description}\\ 
    \hline
    $q$ & $N/m$ & 20000 &  1400 &  Load \\
    $l$ & $m$ & 12 &  0.12 & Length factors for truss members \\
    $A_s$ & $m^2$ & - & $5.9853 \times 10^{-5}$ &  cross-sectional area of steel bars \\
    $A_c$ & $m^2$ & - & 0.0048 &  cross-sectional area of concrete bars \\
    $E_s$ & $Pa$ & $10^{11}$ & $6 \times 10^9$ &  Elastic modulus of steel \\
    $E_c$ & $Pa$ & $2 \times 10^{10}$ & $1.2 \times 10^9$ &  Elastic modulus of concrete\\
    \hline
    \end{tabular}
\end{table}

\autoref{tab:prob6_res} shows the results obtained using different methods. Results corresponding to $\beta_\text{Tar}=2.0$ and $\beta_\text{Tar}=3.0$ are shown.
Compared to the benchmark solution obtained from \cite{Rashki2014}, both TSM and dynamic chaos control (DCC) method are found to yield excellent results.
The other chaos control methods are found to be erroneous.
From the efficiency point of view, TSM is found to be the most efficient with 164 actual function evaluations.

Finally, we compare the advantages gained by each of the modifications (\autoref{tab:p6_diffmod}). 
Since, we have only one probabilistic constraint, {\it Modifications (a) and (c)} have no effect.
{\it Modification (b)} provides significant computational efficiency and without this modification the number of function evaluations increases significantly.

\begin{table}[htbp!]
    \centering
    \caption{Results for the roof truss problem. Results corresponding to DCC, SMCC, hybrid self-adjusted mean value (HSMV) are obtained from \cite{Keshtegar2018}}.
    \label{tab:prob6_res}
    \begin{tabular}{lcccccccccc}
    \hline
    $\beta_\text{Tar}$ & & \multicolumn{3}{c}{2.0} & & \multicolumn{5}{c}{3.0} \\ \cline{3-5}\cline{7-11}
    Method & & DCC & TSM & \cite{Rashki2014} & & SMCC & HSMV & DCC & TSM & \cite{Rashki2014} \\ \hline
    $A_s^* \times 10^-4$ ($m^2$) & & 10.1 & 10.59 & 9.975 & & 6.0 & 6.0 & 11.0 & 11.6 & 10.68 \\
    $A_c^* \times 10^-4$ ($m^2$) & & 3.45 & 3.24 & 3.54 & & 3.0 & 4.19 & 3.87 & 3.57 & 4.05 \\
    $f\left(\bm d \right)$ & & 32.97 & 33.22 & 33.054 & & 23.066 & 27.398 & 36.311 & 36.504 & 36.24 \\
    Opt. iter & & 13 & 11 & - & & 35 & 16 & 12 & 10 & - \\
    $\beta_{\text{MCS}}$ & & 1.989 & 2.02 & 2.02 & & 1.65 & 2.30 & 3.00 & 3.00 & 3.01 \\
    func. calls. & & 1702 & 164 & 40,000 & & 5,224 & 16,437 & 1,102 & 164 & 40,000 \\ \hline
    \end{tabular}
\end{table}

\begin{table}[htbp!]
    \centering
    \caption{Advantage gained using the three modification of TSM (Problem 5). Results corresponding to both $\beta_\text{Tar} = 2.0$ and $\beta_\text{Tar}=3.0$ are shown.}
    \label{tab:p6_diffmod}
    \begin{tabular}{lccc}
    \hline
    \textbf{Methods} & $f\left(\bm d \right)$ & Function calls & Opt. iter \\ \hline \hline
    \textbf{TSM} & 33.22/36.504 & 164 & 11/10 \\ \hline
    TSM (without (a))$^*$ & 33.22/36.504 & 164 & 11/10 \\ 
    TSM (without (b))$^\#$ & 33.05/36.28 & $8 \times  10^6/7 \times  10^6$ & 8/7 \\
    TSM (without (c))$^{\dagger}$ & 133.22/36.504 & 164 & 11/10 \\ \hline
    \multicolumn{4}{l}{$^*$Thresholds for each constraint computed separately.}\\
    \multicolumn{4}{l}{$^\#$No surrogate model is used to approximate the limit state function.} \\
    \multicolumn{4}{l}{$^{\dagger}$Active constraint algorithm has not been utilized.}
    \end{tabular}
\end{table}

\subsection{Crash worthiness of a vehicle due to side impact}\label{subsec:crash_worthiness}
\begin{figure}[htbp!]
    \centering
    \includegraphics[width = 0.4\textwidth]{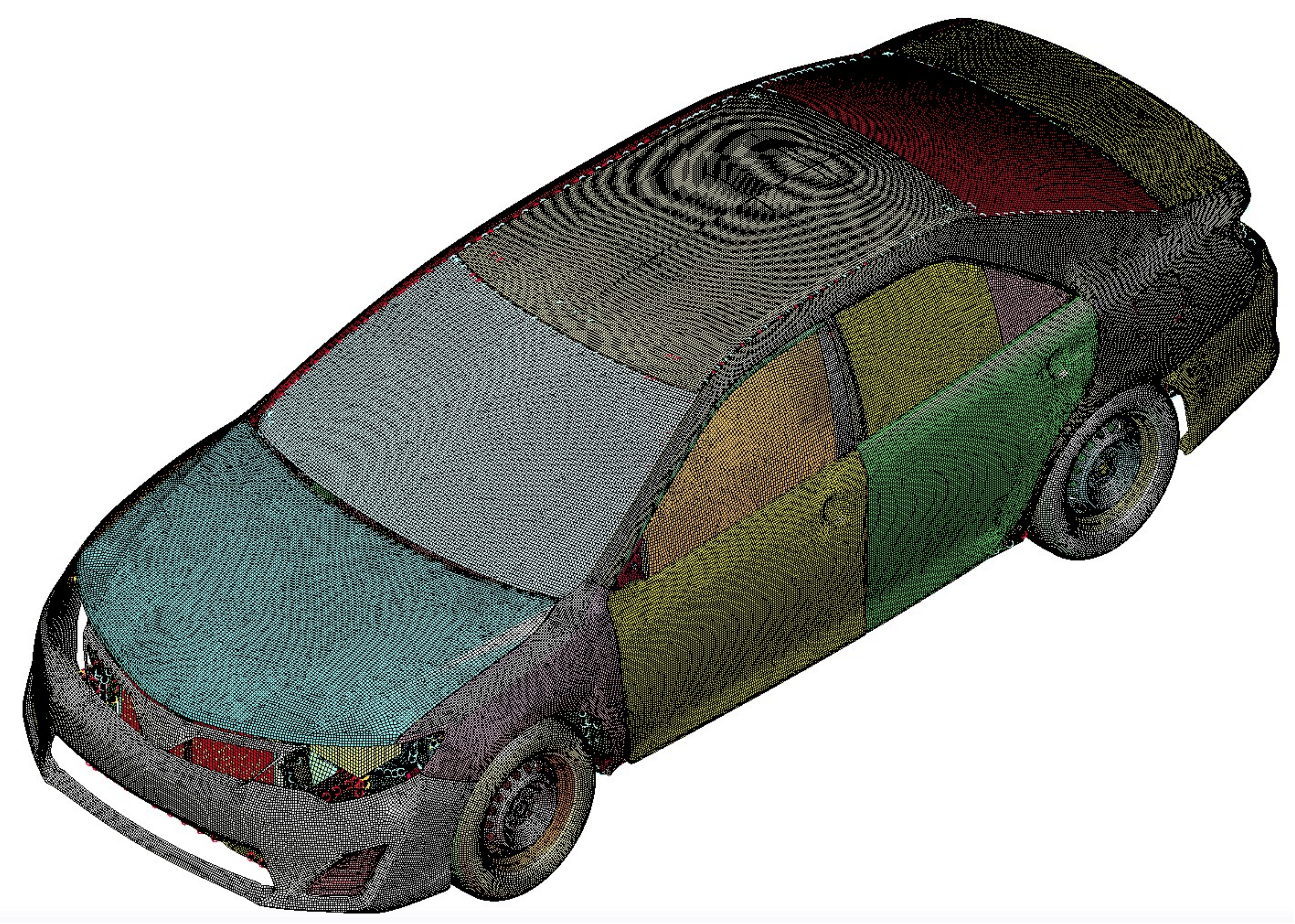}
    \caption{FE model for the vehicle model considered}
    \label{fig:prob2}
\end{figure}

As the last example, we consider RBDO of vehicle against crashworthiness due to side impact.
The system model includes a full-vehicle FE structure
model (\autoref{fig:prob2}), a side impact dummy FE model and side barrier.
The FE model for this problem has been collected from \href{http://www.ncac.gwu.edu/archives/model/index.html}{\color{blue}this link}.
The vehicle model consists of 1.7M nodes and the initial velocity of the barrier is considered to be $49.89$ kmph.
RBDO of this problem has been previously studied in \cite{Youn2004}.
The goal here is to illustrate the applicability of the proposed TSM for solving this problem.

In vehicle design, it is important to satisfy certain side impact requirements.
There exists two primary side impact protection guidelines in literature, (i) National Highway Traffic
Safety Administration (NHTSA) side impact procedures for the
Federal Motor Vehicle Safety Standard (FMVSS) and Canadian
Motor Vehicle Safety Standard (CMVSS) and (ii) European
Enhanced Vehicle-Safety Committee (EEVC) side impact procedure for European vehicle. 
For the present study, we have considered the EEVC side impact criteria.
Details on this criteria can be found in \autoref{tab:crashworthiness}.

\begin{table}[htbp!]
    \centering
    \caption{Safety critera as per EEVC}
    \label{tab:crashworthiness}
    \begin{tabular}{lcc}
    \hline
    \textbf{Performance} &   & \textbf{Top safety rating criteria}  \\ \hline 
    Abdomen load  &  & $\le 1$ \\
     \multirow{3}{*}{Rib deflection} & Upper & \\
      & Middle & $\le 22$ \\
      & Lower & \\
      & & \\
      \multirow{3}{*}{VC (m$/$s)} & Upper & \\
       & Middle & $\le 0.32$ \\
       & Lower & \\
       & & \\
       Public symphysis force (kN) &  & $\le 3$ \\
       HIC &  & $\le 650$ \\ \hline
    \end{tabular}
\end{table}

Based on the safety criteria shown in \autoref{tab:crashworthiness}, the RBDO problem has been formulated as:
\begin{equation}\label{eq:crash_obj}
 \begin{split}
 & \bm d^* = \argmin_{\bm d} w\left(\bm d \right) \\
 & \textit{s.t.}\\ 
 & \mathbb{P} \left( \textit{abdomen load}  \le 1.0 \text{kN}\right) \ge P_s \\
 & \mathbb{P} \left( \text{upper}/\text{middle}/\text{lower} \text{VC}\le 0.32 \text{m}/\text{s} \right) \ge P_s \\
 & \mathbb{P} \left( \text{upper}/\text{middle}/\text{lower} \text{rib deflection}\le 32 \text{mm} \right) \ge P_s \\
 & \mathbb{P} \left(\text{public symphysis force}, F \le 4.0 \text{kN} \right) \ge P_s \\
 & \mathbb{P} \left( \text{velocity of B-pillar at midpoint} \le 9.9 \text{mm}/\text{ms} \right) \ge P_s \\
 & \mathbb{P} \left( \text{velocity of front door at B-pillar}\le 15.7 \text{mm}/\text{ms} \right) \ge P_s \\
 & \bm d^L \le \bm d \le \bm d^U,
 \end{split}
\end{equation}
where $w$ represents the weight of the vehicle.
The system possess eleven random variables, which include thickness, material properties of critical part, barrier height and hitting position. 
Mean of first nine of the eleven variables are considered to be design variables. 
Barrier height and the hitting position are not design variables and may vary from -30 mm to 30 mm
In this problem, all the random variables are assumed to be normally distributed. 
Details of the random variables are provided in \autoref{tab:crash_random}.

\begin{table}[htbp!]
    \centering
    \caption{Statistical parameters of random parameters in \autoref{eq:crash_obj}.}
    \label{tab:crash_random}
    \begin{tabular}{lccccc}
    \hline
    \text{Random variables} & \text{Standrad deviation} & \text{Type} & $\textbf{d}^L$ & $\textbf{d}_i$ & $\textbf{d}^U$\\ 
    \hline
    B-pillar inner (mm) & 0.100 &  Normal & 0.500 &  1.000 & 1.500 \\
    B-pillar reinforce (mm) & 0.100 &  Normal & 0.450 &  1.000 & 1.350 \\
    Floor side inner (mm) & 0.100 &  Normal & 0.500 &  1.000 & 1.500 \\
    Cross member & 0.100 &  Normal & 0.500 &  1.000 & 1.500 \\
    Door beam (mm) & 0.100 &  Normal & 0.875 &  2.000 & 2.625 \\
    Door belt line (mm) & 0.100 &  Normal & 0.400 &  1.000 & 1.200 \\
    Roof rail (mm) & 0.100 &  Normal & 0.400 &  1.000 & 1.200 \\
    Mat. B-pillar inner (GPa) & 0.006 &  Normal & 0.192 &  0.300 & 0.345 \\
    Mat. floor side inner (GPa) & 0.006 &  Normal & 0.192 &  0.300 & 0.345\\
    Barrier height (mm) & 10 &  Normal & \multicolumn{3}{c}{Not design variable}\\
    Barrier hitting (mm) & 10 &  Normal & \multicolumn{3}{c}{Not design variable}\\
    \hline
    \end{tabular}
\end{table}

Using the information provided above, RBDO of the vehicle for crashworthiness due to side-impact has been carried out.
For training the surrogate model ({\it surrogate 1}) in DS-TSM, 256 training samples have been used.
As for {\it surrogate 2}, we started with $32$ training samples.
The other parameters in DS-TSM are considered as specified before.
The RBDO problem has been solved corresponding to target probability of 0.90 and 0.9987.
Since, this is a large-scale problem, generating benchmark solution using brute-force method is difficult.
Therefore, comparison with results in \cite{Youn2004} has been presented.

\autoref{tab:prob2_res} shows the optimum RBDO results corresponding to $P_S = 0.90$, obtained using various methods.
TSM was found to yield the best result (i.e. minimum weight).
The optimum design variables at convergence are shown in \autoref{tab:prob2_res2}. The design variables obtained using various methods are found to be in close proximity.
To further examine the results, the probabilistic constraints at convergence are shown in \autoref{tab:prob2_res3}. 
We observe that all the constraints are satisfied and hence, the design is feasible.

\begin{table}[htbp!]
    \centering
    \caption{Optimized vehicle weight. The RBDO has been carried out corresponding to $P_s = 0.90$. PMA, RIA and MCS results are collected from \cite{Youn2004}.}
    \label{tab:prob2_res}
    \begin{tabular}{ccccc}
       \hline
        \textbf{PMA with HMV} & \textbf{PMA with AMV} & \textbf{RIA with HLRF} & \textbf{MCS} & \textbf{TSM}  \\ \hline
        24.07 & 24.07 & 24.17 & 25.60 & 22.13 \\ \hline
    \end{tabular}
\end{table}

\begin{table}[htbp!]
    \centering
    \caption{Design variables at convergence. The RBDO has been carried out corresponding to $P_s = 0.90$. PMA and RIA results are collected from \cite{Youn2004}.}
    \label{tab:prob2_res2}
    \begin{tabular}{cccccc}
       \hline
        $\bm d $ & \textbf{PMA with HMV} & \textbf{PMA with AMV} & \textbf{RIA with HLRF} & \textbf{MCS} & \textbf{TSM}  \\ \hline
        $d_1$ & 0.5 & 0.5 & 0.5 & - & 0.5 \\
        $d_2$ & 1.309 & 1.309 & 1.309 & - & 1.448 \\
        $d_3$ & 0.5 & 0.5 & 0.5 & - & 0.5 \\
        $d_4$ & 1.240 & 1.239 & 1.239 & - & 0.5 \\
        $d_5$ & 0.610 & 0.610 & 0.667 & - & 0.667 \\
        $d_6$ & 1.5 & 1.5 & 1.5 & - & 1.5 \\
        $d_7$ & 0.5 & 0.5 & 0.5 & - & 0.5 \\
        $d_8$ & 0.345 & 0.345 & 0.345 & - & 0.345 \\
        $d_9$ & 0.192 & 0.192 & 0.258 & - & 0.253 \\ \hline
    \end{tabular}
\end{table}

\begin{table}[htbp!]
    \centering
    \caption{Constraints (as reliability indices) at convergence. The RBDO has been carried out corresponding to $P_s = 0.90$.}
    \label{tab:prob2_res3}
    \begin{tabular}{lcccccccccc}
       \hline
      \textbf{Method}  & $\beta_1$ & $\beta_2$ & $\beta_3$ & $\beta_4$ & $\beta_5$ & $\beta_6$ & $\beta_7$ & $\beta_8$ & $\beta_9$ & $\beta_{10}$ \\ \hline
        {TSM} & 1.85 & $ > 8$ & $ > 8$ & 1.8823 & 2.3871 & 2.8993 & 1.81 & 1.40 & 1.74 & 1.28 \\ \hline
    \end{tabular}
\end{table}

Tables \ref{tab:prob2_res4} and \ref{tab:prob2_res6} shows the optimum RBDO results and the corresponding probabilistic constraints corresponding to $P_s = 0.9987$ obtained using various methods. Although the objective function obtained using the proposed TSM is greater than that obtained using PMA, it is observed that the probabilistic constraints are violated in PMA.
This indicates the PMA yields unsafe solutions.

\begin{table}[htbp!]
    \centering
    \caption{RBDO results corresponding to $P_s = 0.9987$. PMA results are collected from \cite{Youn2004}.}
    \label{tab:prob2_res4}
    \begin{tabular}{lcccccccccc}
    \hline
    \textbf{Methods} & $w\left(\bm d \right)$ & $d_1$ & $d_2$ & $d_3$ & $d_4$ & $d_5$ & $d_6$ & $d_7$ & $d_8$ & $d_9$ \\ \hline
    PMA & 25.39 & 0.511 & 1.417 & 0.5 & 1.352 & 0.658 & 1.473 & 0.5 & 0.345 & 0.192  \\
    TSM  & 26.63 & 1.074 & 1.332 & 0.5 & 0.902 & 1.139 & 1.5 &  0.5 & 0.345 & 0.345  \\
    \hline
    \end{tabular}
\end{table}

\begin{table}[htbp!]
    \centering
    \caption{Constraints (as reliability indices) at convergence. The RBDO has been carried out corresponding to $P_s = 0.9987$.}
    \label{tab:prob2_res6}
    \begin{tabular}{lcccccccccc}
       \hline
      \textbf{Method}  & $\beta_1$ & $\beta_2$ & $\beta_3$ & $\beta_4$ & $\beta_5$ & $\beta_6$ & $\beta_7$ & $\beta_8$ & $\beta_9$ & $\beta_{10}$ \\ \hline
      {PMA} & $>8$ & $ > 8$ & $ > 8$ & 2.12 & 3.35 & 3.1 & 1.60 & $>8$ & 2.60 & 1.32 \\ 
        {TSM} & 3.45 & $ > 8$ & $ > 8$ & 3.20 & 3.14 & $>8$ & 3.00 & 3.62 & 3.02 & 3.57 \\ \hline
    \end{tabular}
\end{table}

The proposed approach requires only 256 actual function calls.
As for convergence, TSM requires only 8 iterations to converge. This is impressive when one considers the fact that the proposed approach involves 11 random variables, 9 design variables and 10 probabilistic constraints.

Finally, we compare the advantages gained by each of the modifications (\autoref{tab:p2_diffmod}). 
Since this problem involves significant computational effort, it is not possible to generate results by excluding {\it Modification (b)}.
It is observed that excluding {\it Modification (a)}, the results obtained are less accurate and the computational efficiency is less.
On the other hand, excluding {\it Modification (c)} doesn't influence the accuracy. 
However, the number of iteration required increases significantly.

\begin{table}[htbp!]
    \centering
    \caption{Advantage gained using the three modification of TSM (Problem 6). Results corresponding to both $P_s = 0.90$ and $P_s=0.9987$ are shown.}
    \label{tab:p2_diffmod}
    \begin{tabular}{lccc}
    \hline
    \textbf{Methods} & $w\left(\bm d \right)$ & Function calls & Opt. iter \\ \hline \hline
    \textbf{TSM} & 22.13/25.14 & 256 & 25/35 \\ \hline
    TSM (without (a))$^*$ & 22.69/26.32 & 256 & 45/60 \\ 
    TSM (without (b))$^\#$ & - & - & - \\
    TSM (without (c))$^{\dagger}$ & 22.13/25.14 & 256 & 85/100 \\ \hline
    \multicolumn{4}{l}{$^*$Thresholds for each constraint computed separately.}\\
    \multicolumn{4}{l}{$^\#$No surrogate model is used to approximate the limit state function.} \\
    \multicolumn{4}{l}{$^{\dagger}$Active constraint algorithm has not been utilized.}
    \end{tabular}
\end{table}

\section{Conclusion}
\label{sec:conclusion}
This paper proposes `threshold shift method' (TSM) for reliability based design optimization (RBDO).
The proposed approach, motivated from sequential optimization and reliability analysis (SORA) method, decouples the RBDO problem. Within each iteration one has to sequentially solve an optimization problem and a reliability analysis problem.
However, unlike SORA, TSM does not utilize a shift-vector; instead, the thresholds of the equivalent constraints are shifted.
It is argued that modifying a constraint, either by shifting the design variables (as in SORA)  or by shifting the threshold of the constraints (as in TSM), influence the other constraints of the system.
Therefore, unlike SORA, the shifted-constraints are determined by solving  a single optimization problem. 
To make TSM scalable, an active-constraint determination scheme has been proposed and integrated with TSM. 
Finally, a practical algorithm for TSM based on two surrogate models has been presented.

The proposed approach has been applied on six benchmark problems from literature. 
Results obtained have been compared with the solutions available in literature. 
It is observed that for all the problems, TSM yields excellent results (converges to optimum solution from very few iterations). 
Moreover, the number of actual function calls (indicative of its computational efficiency) required using TSM is also significantly less -- indicating its possible application to other large scale problems.

\begin{appendices}
\section{PC-Kriging}
\label{appendix-pc_kriging}
 PC-Kriging is novel surrogate model that combines polynomial chaos expansion (PCE) \cite{Pascual2012combined,Xiu2002athe,Sudret2008global,Hu2010adaptive} with Kriging \cite{Bilionis2013multi,Biswas2018saturation,Biswas2016kriging,Mukhopadhyay2016a}.
To be specific, PCE is used to model the mean function within the Kriging surrogate.
This allows one to surrogate a computational model more accurately than PCE or Kriging taken separately.

In PC-Kriging, the trend function, for the $P$-th order ($P>0$) PCE of the function, is mathematically presented as:
\begin{equation}\label{eq:pck}
    \mathcal{M}(\bm{x}) \approx \mathcal{M}_P^{(PCK)}(\bm{x}) = 
    \underbrace {\sum_{\left| \mathbf i = 0 \right|}^P \gamma_{\mathbf i} \psi_{\mathbf i}}_{PCE}
     + \underbrace{{\sigma ^2}Z(\bm{x},\omega )}_{Kriging},
\end{equation}
where $\sigma ^2$ and $Z(\bm{x},\omega)$ are the Gaussian process variance and the zero-mean, unit variance stationary Gaussian process, respectively and $\omega \in \Omega$ denotes an elementary event in the probability space $(\Omega,\mathcal{F},\mathbb{P})$. In \autoref{eq:pck}, $\gamma_{\mathbf i}$ represent the unknown coefficients and $\psi_{\mathbf i}$ represents the orthogonal polynomial of degree $i$. Under limiting conditions, the expression of PC-Kriging in \autoref{eq:pck} converges to either PCE or Kriging.

Despite the superiority of PC-Kriging over PCE and Kriging, it suffers from the {\it curse of dimensionality}.
As the number of input variables $N$ increases, the number of of unknown coefficients to be estimated grows factorially.
This limitation stems from the PCE part of PC-Kriging.
To address this issue, an Optimal PCE-Kriging (OPC-Kriging) was proposed, which utilizes the least angle regression based model selection technique to retain only the important PCE components in the PC-Kriging model.
To be specific, the LAR algorithm results in a list of ranked polynomials which are chosen depending on their correlation to the current residual at each iteration in decreasing order. 
OPC-Kriging consists of an iterative algorithm where each polynomial is added one-by-one to the trend part. 
A flowchart depicting the steps involved in OPC-Kriging is shown in \autoref{fig:opcK}.
For further details on OPC-Kriging, interested readers may refer \cite{Schobi2015polynomial} and \cite{Kersaudy2015a}.
{
\begin{figure}
    \centering
    \includegraphics[width = 0.45\textwidth]{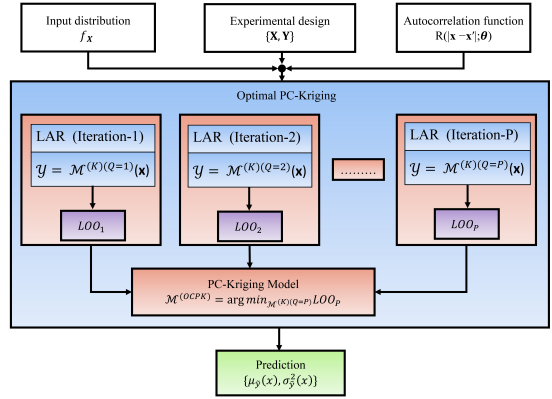}
    \caption{Flowchart depicting the steps involved in OPC-Kriging}
    \label{fig:opcK}
\end{figure}}
\end{appendices}

\section*{Acknowledgement}
SG acknowledges the support of the German Academic Exchange Service (DAAD).

\end{document}